\documentclass[journal ]{new-aiaa}
\usepackage[utf8]{inputenc}
\usepackage{textcomp}

\usepackage{algorithm, algorithmicx, algpseudocode}
\usepackage{bm}
\usepackage{subfig}
\usepackage{amssymb}

\usepackage[svgnames]{xcolor,colortbl}

\usepackage[bookmarks=true,pdftex=false,
bookmarksopen=true,pdfborder={0,0,0},linktocpage=true,
colorlinks=true,linkcolor=RoyalBlue,citecolor=RoyalBlue,
urlcolor=RoyalBlue]{hyperref}

\usepackage{multicol}
\usepackage{multirow}

\usepackage{graphicx}
\usepackage{amsmath}
\usepackage{physics}
\usepackage[version=4]{mhchem}
\usepackage{siunitx}
\usepackage{longtable,tabularx}
\usepackage[capitalize]{cleveref}

\definecolor{rev_color}{rgb}{0.0824, 0.1020, 0.5804}
\newcommand{\hl}{\textcolor{black}}

\DeclareMathOperator*{\argmax}{arg\,max}

\usepackage[acronyms,nonumberlist]{glossaries-extra}
\glsdisablehyper
\setabbreviationstyle[acronym]{long-short}
\renewcommand{\acrfull}[1]{\glsxtrfull{#1}}
\renewcommand{\acrlong}[1]{\glsxtrlong{#1}}
\renewcommand{\Acrshort}[1]{\Glsxtrshort{#1}}
\renewcommand{\acrshort}[1]{\glsxtrshort{#1}}
\renewcommand{\acrshortpl}[1]{\glsxtrshortpl{#1}}
\renewcommand{\acrlongpl}[1]{\glsxtrlongpl{#1}}
\renewcommand{\acrfullpl}[1]{\glsxtrfullpl{#1}}

\newacronym{ads}{ADS}{automatic domain splitting}

\newacronym{cut}{CUT}{conjugate unscented transform}

\newacronym{da}{DA}{differential algebra}

\newacronym{gmm}{GMM}{Gaussian mixture model}

\newacronym{lam}{LAM}{likelihood agreement measure}
\newacronym{loads}{LOADS}{low-order automatic domain splitting}

\newacronym{mc}{MC}{Monte Carlo}

\newacronym{neo}{NEO}{near-Earth object}
\newacronym{nli}{NLI}{nonlinearity index}

\newacronym{pdf}{pdf}{probability density function}
\newacronym{pce}{PCE}{polynomial chaos expansion}

\newacronym{sst}{SST}{state transition tensor}
\newacronym{stm}{STM}{state transition matrix}

\newacronym{ut}{UT}{unscented transform}

\usepackage{comment}
\usepackage{orcidlink}

\newcommand\blfootnote[1]{%
  \begingroup
  \renewcommand\thefootnote{}\footnote{#1}%
  \addtocounter{footnote}{-1}%
  \endgroup
}

\title{A low-order automatic domain splitting approach for nonlinear uncertainty mapping\blfootnote{\textit{This work has been submitted to the Journal of Guidance, Control, and Dynamics for possible publication. Copyright may be transferred without notice, after which this version may no longer be accessible.}}\blfootnote{}}

\author{Matteo Losacco\footnote{Postdoctoral researcher, DCAS,10 Av. Edouard Belin, Toulouse, France (\hl{corresponding author, matteo.losacco0390@gmail.com}.)}\orcidlink{0000-0002-5126-1724}}
\affil{ISAE-SUPAERO, Toulouse, 31400, France}
\author{Alberto Foss\`{a}\footnote{PhD candidate, DCAS,10 Av. Edouard Belin, Toulouse, France.}\orcidlink{0000-0002-0756-4998}}
\affil{ISAE-SUPAERO, Toulouse, 31400, France}
\author{Roberto Armellin\footnote{Professor, Te P\=unaha \=Atea - Space Institute, 20 Symonds Street, Auckland, New Zealand.}\orcidlink{0000-0002-3516-6428}}
\affil{University of Auckland, Auckland, 1010, New Zealand}

\begin{document}

\maketitle

\begin{abstract}
This paper introduces a novel method for the  automatic detection and handling of nonlinearities in a generic transformation. A nonlinearity index that exploits second-order Taylor expansions and  polynomial bounding techniques is first introduced to estimate the Jacobian variation of a nonlinear transformation. This index is then embedded into a \acrlong{loads} algorithm that accurately describes the mapping of an initial uncertainty set through a generic nonlinear transformation by splitting the domain whenever \hl{nonlinearities grow above a predefined threshold}. The algorithm is illustrated in the critical case of orbital uncertainty propagation, and it is coupled with a tailored merging \hl{process} that limits the growth of the domains in time by recombining them when nonlinearities decrease. The \acrlong{loads} algorithm is then combined with \acrlongpl{gmm} to accurately describe the propagation of a \acrlong{pdf}.
\end{abstract}

\section{Introduction}
Orbit uncertainty propagation plays a key role in celestial mechanics. The accurate description of the evolution in time of the uncertainty set of an orbiting object is of crucial importance whenever actions such as follow up observations, data association and orbital conjunctions analyses are required. A rigorous approach to the problem would involve the solution of the Fokker-Plank-Kolmogorov equation to capture the exact evolution in time of the \acrfull{pdf} associated with the orbital state~\citep{Fuller1969}. Although extensive studies have been conducted over the years to identify efficient ways to reduce its computational burden~\citep{Kumar2010,Sun2016}, the approach is still affected by the curse of dimensionality, which makes it unfeasible when propagating uncertainties in a high-dimensional system. As a result, alternative approaches have been developed over the years to reduce the complexity of the problem. According to the level of approximation introduced in the considered dynamical problem, two main categories can be identified, namely nonlinear and linear methods. 

Nonlinear methods can be divided into three groups: sample-based, dynamics-based, and \acrshort{pdf}-based~\citep{Luo2017}. The most common approach among sample-based methods is the classical \acrfull{mc} method. \Acrlong{mc} simulations provide a straightforward way to solve nonlinear and non-Gaussian propagations and reconstruct the propagated \acrshort{pdf}. However, the number of samples required to obtain convergence statistics may be excessively high, thus making the approach computationally heavy. As a result, alternative sampling approaches have been developed over the years, such as the \acrfull{ut}~\citep{Julier2000} and the \acrfull{cut}~\citep{Adurthi2012}. The \acrshort{ut} guarantees a dramatic reduction in the number of required samples but can accurately estimate only the first two moments of the propagated \acrshort{pdf}, while the \acrshort{cut} provides knowledge of higher-order moments but may become \hl{demanding} as the dimension of the problem increases. 

Dynamics-based methods, instead, simplify the uncertainty propagation problem by approximating the transformation. Examples of this kind are the \acrfullpl{sst} method~\citep{Park2012}, the \acrfull{da} method~\citep{Berz1999,Valli2013}, and \acrfull{pce}~\citep{Wiener1938}. Both \acrshortpl{sst} and \acrshort{da} approximate the orbital dynamics by Taylor series expansion, but the latter avoids the derivation \hl{of variational equations and their numerical integration along with the state dynamics}, which makes the approach particularly appealing. Unfortunately, \acrshort{da} offers a local approximation of the mapping of interest, thus failing to provide accurate descriptions for large uncertainty sets or highly nonlinear transformations unless \acrfull{ads} is used \citep{Wittig2015}. Several applications of the combined use of \acrshort{ads} and \acrshort{da} can be found in literature, ranging from the solution of data association problems for optical too-short arcs~\citep{Pirovano2021} to the estimation of the impact probability of \acrlongpl{neo}~\citep{Losacco2018}. The \acrshort{pce} approach, instead, expresses the transformation in the form of a series expansion consisting of orthogonal polynomials, whose coefficients are determined either in a intrusive or non-intrusive way~\citep{Luo2017}.

The last category of nonlinear methods includes \acrshort{pdf}-based methods, e.g. methods based on \acrfullpl{gmm}~\citep{Sorenson1971,Terejanu2008,DeMars2013,Vittaldev2016,Vishwajeet2018}, \hl{which approximate the true state} \acrshort{pdf} \hl{with} a finite sum of weighted Gaussian \acrshortpl{pdf}. The uncertainty propagation problem is thus split into smaller problems for which the linear assumption holds. The technique can be merged with sample-based and dynamics-based methods, obtaining hybrid methods such as the \acrshort{gmm}-\acrshort{sst}~\citep{Fujimoto2015}, the \acrshort{gmm}-\acrshort{da}~\citep{Sun2018}, and the \acrshort{gmm}-\acrshort{pce}~\citep{Vittaldev2016a}. 

The described nonlinear methods generally suffer from either high computational loads or a strong dependency on the problem's dimensionality. As a result, whenever nonlinearities are limited, the so-called linear propagation methods such as the local linearization~\citep{Maybeck1979,Geller2006} and statistical linearization~\citep{Gelb1973,Gelb1974} become appealing. The local linearization technique approximates the considered nonlinear function using the first-order Taylor expansion along the nominal trajectory. The method is simple and computationally efficient, but its accuracy drops in the case of highly nonlinear systems, long term propagations, or significant uncertainties. Similar considerations hold for the statistical linearization, making the range of applicability of linear methods typically limited.

Identifying the boundary between the linear and the nonlinear regimes is not straightforward, as the nature of the problem also depends on the adopted representation of the state. Studies on this topic can be found, for example, in~\citet{Hill2012} and~\citet{Sabol2013}, where the use of equinoctial elements is shown as a valuable tool to help preserve the normality of the propagated \acrshort{pdf}. The availability of an index capable of quantifying the nonlinear nature of the mathematical formulation of the investigated problem can be therefore important, as it could be used as a criterion to select the solving strategy. \citet{Junkins1997} and~\citet{Junkins2004} proposed to use the variation of the Jacobian of the considered nonlinear transformation as a possible \acrfull{nli}. Though representing an appealing way to identify nonlinearities, its validity is strongly affected by the region of interest that is sampled. \citet{Junkins2004} suggest to uniformly sample the worst-case surface, whereas~\citet{Vishwajeet2018} rely on the \acrshort{ut} and sigma points. In both cases, the evaluation of the index requires additional computations with respect to the nominal solution. Conversely, alternative approaches based on the use of Taylor expansions for the determination of Lagrangian coherent structures have been recently proposed~\citep{Perez2015,Tyler2022}.

This paper introduces a new formulation of the \acrshort{nli} to overcome the limitation of sample-based nonlinearity \hl{indices} and exploit the advantages of linear or quasi-linear methods. The index couples DA with polynomial bounding techniques to estimate the Jacobian variation, thus removing any need for sampling. The obtained index is then embedded into a \acrfull{loads} algorithm that accurately describes the mapping of an uncertainty set through a generic nonlinear transform. The technique takes inspiration from the author's high-order \acrlong{ads} algorithm~\citep{Wittig2015}, and reformulates it in a low-order sense by relying on second-order expansions. Starting from an uncertainty set around a nominal solution and a nonlinear transformation, the algorithm exploits the \acrshort{nli} to automatically split the uncertainty whenever nonlinearities are detected. The result of the process is a manifold of quasi-linear representations of the propagated uncertainty set. The algorithm is applied to the case of orbit uncertainty propagation coupled with a novel merging scheme to recombine the polynomials whenever nonlinearities are weaker, thus limiting the total number of domains. 

The natural capability of the \acrshort{loads} algorithm of splitting the uncertainty propagation problem into locally quasi linear problems is then used as a driver for a novel DA-based \acrshort{gmm} method. The resulting \acrshort{loads}-\acrshort{gmm} method is a split/merge technique that automatically tunes the number of \acrshort{gmm} components according to the nonlinearities of the problem detected by the \acrshort{loads} algorithm. As a result, an automatic algorithm for the estimation of the \acrshort{pdf} resulting from the nonlinear transformation of an uncertain variable is obtained. A preview of the combined \acrshort{loads}-\acrshort{gmm} approach has been recently published by the authors for analytical dynamical models~\citep{Fossa2022}. Here we extend the method to a generic framework, and show the application \hl{to} different test cases. \hl{These} test cases are selected to showcase the algorithms' features and performance, rather than to address concrete application scenarios. Analyses  on state representations to reduce nonlinearities are left for future work.

The paper is structured as follows. The \acrshort{da}-based \acrshort{nli} is described in Section~\ref{sec:Index}. Sections~\ref{sec:LOADS} and~\ref{sec:LOADS_unc} offer a general description of the \acrshort{loads} algorithm and then show its implementation for the orbital uncertainty propagation case. The combination of the \acrshort{loads} algorithm with \acrshort{gmm} is then illustrated in Section~\ref{sec:LOADS-GMM}, and the major modifications with respect to the standard \acrshort{loads} are described. Finally, the proposed methods are applied to \hl{a number of} test cases in Section~\ref{sec:Sim}, where the performance of the algorithm and its sensitivity to \hl{the} control parameters are studied.

\section{The \acrlong{da}-based \acrlong{nli}}
\label{sec:Index}
Consider a generic multivariate nonlinear transformation \hl{written in non-dimensional quantities}
\begin{equation}
\bm{y} = \bm{f}(\bm{x})
\end{equation}
\hl{and two multivariate random variables $\bm{X}$ and $\bm{Y}$, defined on} $\mathbb{R}^n$ \hl{and} $\mathbb{R}^m$ \hl{and with \acrshort{pdf} $p_{\bm{X}}$ and $p_{\bm{Y}}$, respectively}, while $\bm{f}$ is non-odd and $\mathcal{C}^2$ over $\bm{x}$. In the frame of \acrshort{da}, \hl{$\bm{X}$} can be expressed as a function of both its \hl{expected} value and the associated uncertainty. As a result, the following expression can be obtained
\begin{equation}
\label{eq:x}
[x_j] = \bar{x}_j+\beta_j\delta x_j\qquad j=1,\ldots,n
\end{equation}
Equation~\eqref{eq:x} indicates that the \acrshort{da} quantity $[x_j]$, conventionally reported within square brackets, is the sum of two terms. The first term $\bar{x}_j$ represents the \hl{expected or} nominal value of the component $j$, whereas the second term is the product of two quantities, namely $\beta_j$ and $\delta x_j$. The coefficient $\beta_j\hl{\geq 0}$ is a constant corresponding to the semi-amplitude of the confidence interval of the component $j$. This confidence interval is defined as the region around the nominal value that is expected to include a prescribed level of probability mass, and is an input of the problem. The term $\delta x_j$ is the small variation around $\bar{x}_j$, so that $\delta x_j\in[-1,+1]$. The expression can be written as
\begin{equation}
[x_j] = \mathcal{T}^{(1)}_{x_j}(\delta x_j)
\end{equation}
where $\mathcal{T}^{(1)}_{x_j}$ indicates the first-order (1 superscript) Taylor expansion of $x_j$ (subscript), expressed as a function of the deviation $\delta x_j$. Let us now apply the nonlinear transformation $\bm{f}$. In the \acrshort{da} framework, we have that $[\bm{y}]=\bm{f}([\bm{x}])$. If we restrict the analysis to second-order expansions, each component $[y_i]$\hl{, with $i=1,\ldots,m$,} can be expressed as

\begin{equation}
\label{eq:y_i}
\begin{aligned}
\relax[y_i] = a_i + a_{i,1}\delta x_1+\ldots+a_{i,n}\delta x_n+a_{i,11}\delta x_1^2+a_{i,12}\delta x_1\delta x_2+a_{i,13}\delta x_1\delta x_3+\ldots+a_{i,1n}\delta x_1\delta x_n+\\
+a_{i,22}\delta x_2^2+a_{i,23}\delta x_2\delta x_3+\ldots+a_{i,2n}\delta x_2\delta x_n+\\
\vdots\ \\
+a_{i,nn}\delta x_n^2\
\end{aligned}
\end{equation}
\hl{or, more compactly
\begin{equation}
\begin{aligned}
    \relax[y_i] &= a_i+\sum\limits_{p=1}^n a_{i,p}\delta x_p+\sum\limits_{p=1}^n\sum\limits_{k=p}^n a_{i,pk}\delta x_p\delta x_k\\
    &=\mathcal{T}^{(2)}_{y_i}(\delta x_1,\ldots,\delta x_n) = \mathcal{T}^{(2)}_{y_i}(\delta\bm{x})
\end{aligned}
\qquad i=1,\ldots,m
\end{equation}
}
That is, $[y_i]$ can be written as a second-order Taylor expansion $\mathcal{T}^{(2)}_{y_i}$ in $\delta\bm{x}$. Let us now consider the Jacobian of the transformation. The generic term $[J_{ij}]$ reads
\begin{equation}
[J_{ij}]=\dfrac{\partial[y_i]}{\beta_j\partial\delta x_j}
\end{equation}
where the term $\beta_j$ has been introduced to remove the dependency on the initial scaling. By replacing $[y_i]$ with its Taylor series expansion given by Eq.~\eqref{eq:y_i}, we can write
\hl{
\begin{equation}
\begin{aligned}
    \relax[J_{ij}] &= \dfrac{1}{\beta_j}\left(a_{i,j}+a_{i,1j}\delta x_1+\ldots+a_{i,j-1j}\delta x_{j-1}+2a_{i,jj}\delta x_j+a_{i,jj+1}\delta x_{j+1}+\ldots+a_{i,jn}\delta x_n\right)\\
    &=\dfrac{a_{i,j}}{\beta_j}+\sum\limits_{p=1}^{j-1}\dfrac{a_{i,pj}}{\beta_j}\delta x_p+2\dfrac{a_{i,jj}}{\beta_j}\delta x_j+\sum_{p=j+1}^n \dfrac{a_{i,jp}}{\beta_j}\delta x_p
    =\bar{J}_{ij}+\sum\limits_{p=1}^n c_{ij,p}\delta x_p
\end{aligned}
\label{eq:J_ij}
\end{equation}
}
The single term $[J_{ij}]$ is made of a constant part $\bar{J}_{ij}=a_{i,j}/\beta_j$ plus a first-order Taylor expansion in the deviations $\delta \bm{x}$. Equation~\eqref{eq:J_ij} can be therefore written as
\begin{equation}
\label{eq:J_ij2}
\left[J_{ij}\right] = \bar{J}_{ij}+\delta J_{ij}(\delta\bm{x})=\bar{J}_{ij}+\mathcal{T}^{(1)}_{\delta J_{ij}}(\delta\bm{x})=\mathcal{T}^{(1)}_{J_{ij}}(\delta\bm{x})
\end{equation}
or \hl{in matrix form as}
\hl{
\begin{equation}
    [\bm{J}] = \bar{\bm{J}}+\delta\bm{J}(\delta\bm{x})=\bar{\bm{J}}+\mathcal{T}_{\delta\bm{J}}^{(1)}(\delta\bm{x})=\mathcal{T}^{(1)}_{\bm{J}}(\delta\bm{x})
\end{equation}
}
%
%
%
Like its individual terms, the matrix $[\bm{J}]$ is the sum of a constant part $\bm{\bar{J}}$ plus a first-order term $\delta\bm{J}$. The constant term $\bar{\bm{J}}$ is the Jacobian associated with the nominal solution $\bar{\bm{x}}$, and can be computed by considering null deviations, i.e.
\begin{equation}
\bar{\bm{J}}=\mathcal{T}^{(1)}_{\bm{J}}(\delta\bm{x}=\bm{0})
\end{equation}
Conversely, the Jacobian $\bm{J}^*$ associated with a deviation $\delta \bm{x}^*$ can be retrieved by mapping this deviation with $\mathcal{T}_{\bm{J}}(\delta\bm{x})$, i.e.
\begin{equation}
\bm{J}^*=\mathcal{T}_{\bm{J}}(\delta\bm{x}^*)
\end{equation}
Let us now consider the non-constant part of the Jacobian matrix, $\delta\bm{J}$. In the case of a linear transformation, the Jacobian \hl{is} constant, thus all terms $\delta J_{ij}$ \hl{are} equal to zero. For a generic transformation, this is no longer true, thus each term $[J_{ij}]$ includes both a constant part plus a $\delta J_{ij}$ element. \hl{
If we assume that the nonlinear part of the Taylor expansion of $\bm{f}$ is dominated by second-order terms, then the magnitude of the $\delta J_{ij}$ terms is expected to increase with the nonlinearity of the transformation, i.e. the $c_{ij,p}$ coefficients grow as the linear assumption becomes less and less valid. As a result, if} we succeed in quantifying the overall magnitude of all the $\delta J_{ij}$ terms, then an automatic measure of the nonlinear nature of the problem can be obtained. But how to compute this magnitude? The answer is given by polynomial bounding techniques. More specifically, if we recall the expression for $\delta J_{ij}$, i.e.
\begin{equation}
\label{eq:dJ_ij}
\delta J_{ij} = \sum\limits_{p=1}^n c_{ij,p}\delta x_p
\end{equation}
we can see that each \hl{entry} $\delta J_{ij}$ is the weighted sum of $n$ \hl{monomials $c_{ij,p}\delta x_p$}. By assumption, \hl{the terms $\delta x_p$} can vary between -1 and +1. As a result, an upper bound $b_{\delta J_{ij}}$ for $\delta J_{ij}$ can be obtained as the sum of the absolute values of the $c_{ij,p}$ \hl{coefficients}, i.e.
\begin{equation}
\label{eq:B_ij}
b_{\delta J_{ij}} = \sum_{p=1}^n \left\lvert c_{ij,p}\right\rvert
\end{equation}
while the lower bound is equal to $-b_{\delta J_{ij}}$. This means that the range of $\delta J_{ij}$ is simply a symmetric interval centred in 0 with semi-amplitude $b_{\delta J_{ij}}$, i.e. $\delta J_{ij}\in[-b_{\delta J_{ij}},b_{\delta J_{ij}}]$. By performing the same computation for all $\delta J_{ij}$ terms, a matrix of upper bounds $\bm{B}_{\delta\bm{J}}$ can be built as
\begin{equation}
\bm{B}_{\delta\bm{J}} = \begin{bmatrix}
b_{\delta J_{11}}& \cdots& b_{\delta J_{1n}}\\
\vdots& \vdots& \vdots\\
b_{\delta J_{m1}}& \cdots& b_{\delta J_{mn}}
\end{bmatrix}
\end{equation}
A possible measure of $\bm{B}_{\delta\bm{J}}$ can be obtained by computing its Frobenius norm, i.e.
\begin{equation}
\label{eq:B_norm}
\left\lVert \bm{B}_{\delta\bm{J}} \right\rVert_2=\sqrt{\displaystyle\sum_{i=1}^m\sum_{j=1}^n b^2_{\delta J_{ij}}}
\end{equation}
The $\left\lVert \bm{B}_{\delta\bm{J}} \right\rVert_2$ norm sums up the squares of the $b_{\delta J_{ij}}$ terms. In case of a linear transformation, these terms are zero, thus $\left\lVert \bm{B}_{\delta\bm{J}} \right\rVert_2=0$. If the transformation is nonlinear, \hl{under the described assumptions,} the $b_{\delta J_{ij}}$ terms increase with an increase in the nonlinearity of the transformation, thus $\left\lVert \bm{B}_{\delta\bm{J}} \right\rVert_2$ increases as well. If we normalize  this measure by the Frobenius norm of the constant part of the Jacobian, the following ratio can be computed
\begin{equation}
\label{eq:nu}
\nu = \dfrac{\left\lVert \bm{B}_{\delta\bm{J}} \right\rVert_2}{\left\lVert \bar{\bm{J}}\right\rVert_2}=\sqrt{\dfrac{\displaystyle\sum_{i=1}^m\sum_{j=1}^n b^2_{\delta J_{ij}}}{\displaystyle\sum_{i=1}^m\sum_{j=1}^n  \bar{J}_{ij}^2}} 
\end{equation}
By construction, this ratio increases with the nonlinearity of the considered transformation \hl{over the considered set}, thus can be seen as a nondimensional measure of nonlinearity. From now on, it will be referred to as \acrfull{nli} or $\nu$.
As a result, given a generic transformation $\bm{f}$ and \hl{an uncertainty set} $[\bm{x}]$, an estimate of the nonlinear nature of the mathematical formulation of the problem \hl{over this set} can be obtained by inspecting the value of $\nu$. More specifically, by imposing a non-dimensional nonlinearity threshold $\varepsilon_{\nu}$, we can say that
\begin{equation}
\begin{cases}
\nu\leq\varepsilon_{\nu}\rightarrow \textrm{\hl{quasi-}linear}\\
\nu>\varepsilon_{\nu}\rightarrow \textrm{nonlinear}
\end{cases}
\end{equation}
The expression of the \acrshort{nli} here described closely resembles the formulation proposed by~\hl{\citet{Junkins2004}}, and can be seen as its reformulation in \acrshort{da} sense. The proposed approach has however a clear advantage. By relying on \acrshort{da} and using the polynomial bounder, an estimate of the \acrshort{nli} can be obtained rigorously without relying on any sampling technique. A more detailed \hl{comparison} will be presented in Section~\ref{sec:Sim}.

\section{The LOADS algorithm}
\label{sec:LOADS}
Consider the nonlinear transformation $\bm{f}$ and the uncertainty  set $[\bm{x}]$ and assume that $\nu>\varepsilon_{\nu}$. Intuitively, a better agreement with a linear behaviour is obtained by recomputing the index with the same function but on smaller sets. \hl{Let us consider Eq.~\eqref{eq:J_ij} and assess the effect of an uncertainty set size change. More specifically, let us consider a transformation $\delta\bm{x} = \bm{k}\odot\delta\tilde{\bm{x}}$, where $\odot$ is the Hadamard or component-wise product, $\bm{k}\in\mathbb{R}^n_{>0}$ and $\delta\tilde{\bm{x}}\in[-1,1]^n$. The transformation essentially stretches or shrinks the amplitude of the uncertainty set around the nominal solution, which is now equal to $\beta_jk_j$ for each component $j$. If we compose Eq.~\eqref{eq:J_ij} with $\delta\bm{x}$, we obtain} 
\hl{\begin{equation}
\begin{aligned}
    \left[\tilde{J}_{ij}\right] &= \dfrac{a_{i,j}}{\beta_j}+\sum\limits_{p=1}^{j-1}\dfrac{a_{i,pj}}{\beta_j}k_p\delta \tilde{x}_p+2\dfrac{a_{i,jj}}{\beta_j}k_j\delta \tilde{x}_j+\sum_{p=j+1}^n \dfrac{a_{i,jp}}{\beta_j}k_p\delta \tilde{x}_p\\
    &=\bar{J}_{ij}+\sum\limits_{p=1}^n \tilde{c}_{ij,p}\delta \tilde{x}_p
\end{aligned}
\label{eq:J_ij_comp}
\end{equation}}
\hl{
From the equation above it is clear that if the size of the uncertainty set is progressively reduced (i.e., when $\bm{k}$ gets smaller) the magnitude of the non-constant part of the Jacobian becomes smaller, and the assumption of locally linear behaviour becomes more and more valid.}

Starting from this consideration, a splitting algorithm is introduced. The method, defined \acrfull{loads}, represents a reformulation in low-order sense of the automatic domain splitting (ADS) method proposed by the authors in~\citet{Wittig2015} for higher orders. The \acrshort{loads} method accurately maps an uncertainty set through a nonlinear transformation by relying on a sequence of second-order expansions. The introduced \acrshort{nli} $\nu$ represents the kernel of the method. When this index exceeds an imposed threshold, the linear assumption is no longer met. Thus a split is performed in a determined direction.
The splitting direction is estimated by identifying the variable $\delta x_e$ that contributes the most to the nonlinearity. More specifically, let us define with $\delta\bm{x}_e$ the $n$-dimensional deviation vector whose only non-zero component is $\delta x_e\hl{\in[-1,1]}$, that is
\begin{equation}
\label{eq:delta}
\delta\bm{x}_e = \left\{0,\ldots,0,\delta x_e,0,\ldots,0\right\}^{\textrm{T}}\qquad e\in[1,\ldots,n]
\end{equation}
Given the expression of the generic element $[J_{ij}]$ of Eq.~\eqref{eq:J_ij},  the composition with $\delta\bm{x}_e$ results into
\begin{equation}
\label{eq:J_j}
[\left.J_{ij}\right\rvert_e]=[J_{ij}]\circ \delta\bm{x}_e = \bar{J}_{ij}+c_{ij,e}\delta x_e
\end{equation}
As can be seen, the constant part of the $[J_{ij}]$ term is unaltered by the composition, whereas the non-constant part now contains only the deviation in the $\delta x_e$ variable. As a result, the associated directional \acrshort{nli} $\nu_e$ is computed as
\begin{equation}
\label{eq:nu_e}
\nu_e = \sqrt{\cfrac{\displaystyle\sum\limits_{i=1}^m \sum\limits_{j=1}^n b^2_{\delta\left.J_{ij}\right\rvert_e}}{\displaystyle\sum\limits_{i=1}^m\sum\limits_{j=1}^n \bar{J}_{ij}^2}}=\sqrt{\cfrac{\displaystyle\sum\limits_{i=1}^m\sum\limits_{j=1}^n c_{ij,e}^2}{\displaystyle\sum\limits_{i=1}^m\sum\limits_{j=1}^n \left(\dfrac{\partial [y_i]}{\beta_j\partial\delta x_j }(\bm{0})\right)^2}}
\end{equation}
The splitting direction $d$ is then selected as the one corresponding to the maximum directional index, i.e.
\begin{equation}
\label{eq:dir}
    d = \argmax_e\{\nu_e\}
\end{equation}
Once the splitting direction is defined, the initial set $[\bm{x}]$ is split into three subsets. More specifically, given
\begin{equation}
\begin{aligned}
\label{eq:split_1}
\delta\bm{x}_{d}^{(s+1)} &= \left\{\hl{\delta x_1,\ldots,\delta x_{d-1}},\dfrac{1}{3}\delta x_d+\dfrac{2}{3}\left(s-1\right),\hl{\delta x_{d+1},\ldots,\delta x_n}\right\}^{\textrm{T}}\qquad s=\{0,1,2\}
\end{aligned}
\end{equation}
three new sets are obtained as
\begin{equation}
\begin{aligned}
\label{eq:split_2}
\left[\bm{x}^{(s+1)}\right]&= \left[\bm{x}\right]\circ\delta\bm{x}_{d}^{(s+1)}\qquad \beta_{d}^{(s+1)} = \dfrac{1}{3}\beta_d
\end{aligned}
\end{equation}

\begin{algorithm}[!t]
\begin{algorithmic}
\setstretch{1.2}
    \Function{loads}{$\bm{f},[\bm{x}],\varepsilon_\nu,N_{max,e}\forall e$}
        \State Initialize $W_{\bm{x}}$, $M_{\bm{x}}, M_{\bm{y}}$
        \State Add $[\bm{x}]$ to $W_{\bm{x}}$
        \While{$W_{\bm{x}}$ is not empty}
            \State Remove the first set $[\bm{x}^{(k)}]$ from $W_{\bm{x}}$
            \State Evaluate $[\bm{y}^{(k)}]$ \Comment{see \cref{eq:y_i}}.
            \State Compute the \acrshort{nli} $v^{(k)}$ \Comment{see \cref{eq:J_ij,eq:dJ_ij,eq:B_ij,eq:B_norm,eq:nu}}
            \If{$\nu^{(k)}\leq\varepsilon_{\nu}$}
                \State Add $[\bm{x}^{(k)}]$ to $M_{\bm{x}}$ and $[\bm{y}^{(k)}]$ to $M_{\bm{y}}$
            \Else
                \State Compute the splitting direction $d$ \Comment{see \cref{eq:delta,eq:J_j,eq:nu_e,eq:dir}}
                \If{$N_d^{(k)}<N_{max,d}$}
                    \State Split $[\bm{x}^{(k)}]$ into $[\bm{x}^{(k+s)}]$ for $s=0,1,2$ \Comment{see \cref{eq:split_1,eq:split_2}}
                    \State Add all $[\bm{x}^{(k+s)}]$ to $W_{\bm{x}}$
                \Else
                    \State Add $[\bm{x}^{(k)}]$ to $M_{\bm{x}}$ and $[\bm{y}^{(k)}]$ to $M_{\bm{y}}$
                \EndIf
            \EndIf
        \EndWhile
        \State \Return $M_{\bm{x}},M_{\bm{y}}$
    \EndFunction
\end{algorithmic}
\caption{\glsentryshort{loads} algorithm}
\label{alg:loads_algorithm}
\end{algorithm}

\hl{The process starts by initializing a so-called working manifold $W_{\bm{x}}$ as $W_{\bm{x}} = \left\{[\bm{x}]\right\}$. The algorithm then extract $[\bm{x}]$, computes its NLI, and performs a check. If $\nu$ is lower than the imposed threshold $\varepsilon_{\nu}$, the requirements are satisfied, and the algorithm stops. Otherwise, the set is split, and the three generated sets are stored into $W_{\bm{x}}$}. The procedure is then repeated on the generated sets, and it terminates when each set has $\nu<\varepsilon_{\nu}$ or reaches a minimum box size. \hl{This second condition is enforced by limiting the number of splits $N_d$ along direction $d$ to a predefined maximum number $N_{max,d}$.}
A summary of the method is given in Algorithm~\ref{alg:loads_algorithm}. At the end of the process, the algorithm \hl{returns two manifolds}. The first \hl{one}, here referred as $M_{\bm{x}}$, includes all the subsets the initial $[\bm{x}]$ is divided into, i.e.
\begin{equation}
M_{\bm{x}} = \left\{\left[\bm{x}^{(k)}\right]:\bigcup\limits_{k=1}^{n_{set}}\left[\bm{x}^{(k)}\right] = \left[\bm{x}\right]\right\}
\end{equation}
where $n_{set}$ is the overall number of sets. The notation $\left[\bm{x}^{(k)}\right]$ indicates the generic set $k$ with splitting history $H^{(k)}$. This splitting history can be expressed as a series of couples, whose terms indicate the splitting direction and the splitting side (``0'' for left, ``1'' for middle, or ``2'' for right), \hl{respectively}. As a result, if we  indicate with $N^{(k)}$ the overall number of splits for set $k$, $H^{(k)}$ can be expressed as 
\begin{equation}
    H^{(k)}=\left(d^{(k)}_1 s^{(k)}_1;\ldots;d^{(k)}_{N^{(k)}}s^{(k)}_{N^{(k)}}\right)
\end{equation}
where $d^{(k)}_1$ is the splitting direction for set $k$ at the first splitting event, $s^{(k)}_1$ is the corresponding side, etc. In the rest of the paper, the notations $\left[\bm{x}^{(k)}\right]$ (index notation) and $\left[\bm{x}^{\left(H^{(k)}\right)}\right]$ (splitting history notation) will be used interchangeably to indicate the generic element $k$ of a manifold.

The second manifold, here referred to as $M_{\bm{y}}$, collects all the images of the generated $\left[\bm{x}^{(k)}\right]$ sets, i.e.
\begin{equation}
M_{\bm{y}} = \left\{\left[\bm{y}^{(k)}\right]: \left[\bm{y}^{(k)}\right]=\bm{f}\left(\left[\bm{x}^{(k)}\right]\right)\right\}\qquad k=1,\ldots,n_{set}
\end{equation}

\hl{As previously described, the parameters $\varepsilon_{\nu}$ and $N_{max,d}$ play a key role as they govern when splits are required, and how many splits can be performed. This aspect is crucial, as the number of sets may significantly grow with the dimensionality of the problem. The two parameters, therefore, shall be the result of a trade-off between desired accuracy and required computational load. An analysis of the role of $\varepsilon_{\nu}$ on the \acrshort{loads} performance is shown in Section~\ref{subsub:LOADS_param}.}

\hl{A final remark shall be made on the formulation of the \acrshort{nli}. As described in Section~\ref{sec:Index}, its derivation relies on the assumption that second-order terms dominate the nonlinear portion of the Taylor expansion of $\bm{f}$. This holds for all the transformations illustrated in this paper and is typical within the convergence radius of convergent power series, where the size of the coefficients follows a trend of exponential decay~\citep{Wittig2015}. In this circumstance, an increase in nonlinearity causes the $c_{ij,p}$ terms to grow, and so $\nu$ quantifies the nonlinearity of $\bm{f}$. In the general case, however, this may not be true, as second-order terms may not dominate. In such a scenario our \acrshort{nli} can detect, but not measure, the nonlinearites in $\bm{f}$, thus providing just an insight on their level. Nevertheless, based on this estimate, the \acrshort{loads} algorithm can start splitting the uncertainty set, thus scaling down the expansion coefficients by a factor $k^n$, where $k$ is the scaling factor and $n$ is the order each coefficient refers to. That is, coefficients become smaller and smaller as the set shrinks, and the higher the order, the larger its reduction factor. As a result, as splits proceed, the Taylor expansion is progressively reduced to a convergent series for which second-order terms are now dominant and $\nu$ is an accurate measure of nonlinearity.}
%
 %
\hl{Overall, though \acrshort{nli} may not be an exact measure of nonlinearity over the whole domain for any generic nonlinear transformation $\bm{f}$, it can always detect nonlinearities and therefore it can be used as a tool to generate subsets over which the linearity assumption locally holds true.}

\section{Orbital uncertainty propagation}
\label{sec:LOADS_unc}
The \acrshort{loads} algorithm is valid for any regular nonlinear transformation. This section illustrates its implementation in the critical case of orbital uncertainty propagation. More specifically, consider the dynamics
\begin{equation}
\dot{\bm{x}}(t) = \bm{g}(\bm{x}(t),t)\qquad \bm{x}(t_0)=\bm{x}_0
\label{eq:ivp_orbit_up}
\end{equation}
Let us also define for convenience $\bm{G}(\bm{x},t)=\int\limits_{t_0}^{t_f} \bm{g}(\bm{x},t)\textrm{d}t$, with $\bm{x}(t_0)=\bm{x}_0$ such that $\bm{G}(\bm{x},t_0)=\bm{x}_0$ and $\bm{G}(\bm{x}_0,t_f)=\bm{x}_f$. Without loss of generality, let us assume that $\hl{\bm{X}_0}\sim\mathcal{N}\left(\bm{\mu}_0,\bm{P}_{0}\right)$, i.e. the initial state is normally distributed with mean $\bm{\mu}_0$ and covariance $\bm{P}_0$. By defining with $\bm{V}$ and $\bm{\Lambda}=\textrm{diag}\{\lambda_j\}$ the matrices of eigenvectors and eigenvalues of $\bm{P}_0$, respectively, the set of initial conditions can be reformulated \hl{in the DA framework} as
\begin{equation}
[\bm{x}_0] = [\bm{x}(t_0)]  = \bm{\mu}_0+\bm{V}
\begin{Bmatrix}
\alpha_{\lambda_1}\sqrt{\lambda_1}\delta x_1&
\ldots&
\alpha_{\lambda_n}\sqrt{\lambda_n}\delta x_n
\end{Bmatrix}^{\textrm{T}}
\end{equation}
where $\alpha_{\lambda_j}\sqrt{\lambda_j}$ is a scaling factor determining the size of the uncertainty set along direction $j$. The $\alpha_{\lambda_j}$ is a nonnegative factor governing the amount of probability mass included in direction $j$, and is typically selected as \hl{an integer value, e.g. 2, 3, or 4}. Let us indicate with $\varepsilon_{\nu}$ the imposed nonlinearity threshold, \hl{with} $N_{max,e}$ the maximum number of splits for the direction $e$, and \hl{with} $t_f$ the final integration time. Our goal is to propagate the uncertainty set $[\bm{x}_0]$ from $t_0$ to $t_f$ by exploiting the \acrshort{loads} algorithm with the imposed constraints on accuracy and maximum number of splits. In this work, \cref{eq:ivp_orbit_up} is integrated numerically using a \acrshort{da}-aware implementation of the 8(5,3) Dormand-Prince integrator\footnote{\href{https://www.hipparchus.org/apidocs/org/hipparchus/ode/nonstiff/DormandPrince853FieldIntegrator.html}{https://www.hipparchus.org/apidocs/org/hipparchus/ode/nonstiff/DormandPrince853FieldIntegrator.html}}.

Let us then define the working manifold $W_{t}$
\begin{equation}
W_t = \left\{\left[\bm{x}^{(k)}(t_k)\right]:\left( N_e^{(k)}<N_{max,e}\,\forall e=1,\ldots,n\right)\land\left(t_k<t_f\right)\right\}
\end{equation}
The manifold contains all the $k$ sets propagated at time $t_k$ that have a number of splits per direction lower than the maximum allowed ($N_e^{(k)}<N_{max,e}\,\forall e=1,\ldots,n$), and have not \hl{yet reached} the final propagation time ($t_k<t_f$).
At the very beginning of the simulations, only one set is \hl{present}, thus $W_{t}=\left\{\left[\bm{x}_0\right]\right\}$.

The \acrshort{loads} will then generate two new manifolds, collecting the final polynomial maps at $t_f$ and the corresponding subsets in the initial domain, respectively. The first manifold is called $M_{f}$ and can be expressed as
\begin{equation}
M_f=\left\{\left[\bm{x}^{(k)}(t_k)\right]: t_k= t_f\right\}\qquad k=1\ldots,n_{set}
\end{equation}
The manifold $M_{f}$ collects all the sets propagated to the final epoch $t_f$. These include both sets that satisfy the nonlinearity constraint and domains that reach the maximum number of split along one direction at $t<t_f$ and are then propagated to $t_f$ with no more splits allowed.

The second manifold is called $M_0$ and collects all the preimages of the elements of $M_f$, i.e.
\begin{equation}
M_{0}=\left\{\left[\bm{x}^{(k)}(t_0)\right]: \left[\bm{x}^{(k)}(t_f)\right]=\bm{G}\left(\left[\bm{x}^{(k)}(t_0)\right],t_f\right)\forall\left[\bm{x}^{(k)}(t_f)\right]\in M_f\right\}
\end{equation}
where $\bigcup_{k}\left[\bm{x}^{(k)}(t_0)\right]=\left[\bm{x}_0\right]$. The manifold is automatically retrieved from $[\bm{x}_0]$ and $M_f$ \hl{at the end of} the propagation.

The iterative procedure is detailed in Algorithm \ref{alg:loads_algorithm_up}. At each integration step, a \hl{check} on the \hl{\acrshort{nli} of the current set} is made, and if \hl{it exceeds} the prescribed \hl{threshold}, \hl{a split} is performed. \hl{The index} $\nu$ is built on the Taylor expansion of the \acrfull{stm}, which now plays the role of the Jacobian of the transformation. As an example, consider the generic set $\left[\bm{x}^{(k)}(t_k)\right]$ of $W_t$ as available at time epoch $t_k$, and then propagate it to $t_k+\Delta t$, thus obtaining $\left[\bm{x}^{(k)}(t_k+\Delta t)\right]$. The DA expansion of the \acrshort{stm} can be obtained as
\begin{equation}
\label{eq:STM}
\left[\bm{J}^{(k)}\right] = \left[\dfrac{\partial G_i\left(\left[\bm{x}^{(k)}\left(t_k\right)\right],t_k+\Delta t\right)}{\alpha_{\lambda_j}\sqrt{\lambda^{(k)}_j}\partial \delta x_j}\right]=\left[\dfrac{\partial \left[x^{(k)}_i(t_k+\Delta t)\right]}{\alpha_{\lambda_j}\sqrt{\lambda^{(k)}_j}\partial \delta x_j}\right]=\left[\bm{\Phi}^{(k)}(t_k+\Delta t,t_0)\right]
\end{equation}
The \hl{\acrlong{nli} is then built from \cref{eq:STM} as follows. First, the non-constant part of the \acrshort{stm} is retrieved as}
\begin{equation}
\label{eq:dSTM}
\delta\bm{\Phi}^{(k)}(t_k+\Delta t,t_0)=\left[\bm{\Phi}^{(k)}(t_k+\Delta t,t_0)\right]-\bar{\bm{\Phi}}^{(k)}(t_k+\Delta t,t_0)=\left[\displaystyle\sum\limits_{p=1}^n c_{ij,p}^{(k)}\delta x_p\right]
\end{equation}
\hl{As a result, the} matrix of upper bounds $\bm{B}_{\delta\bm{\Phi}^{(k)}(t_k+\Delta t,t_0)}$ \hl{is readily obtained as}
\begin{equation}
\label{eq:dSTM_bounds}
\bm{B}_{\delta\bm{\Phi}^{(k)}(t_k+\Delta t,t_0)}=\left[ \displaystyle\sum\limits_{p=1}^n \left\lvert c_{ij,p}^{(k)}\right\rvert\right]
\end{equation}
\begin{algorithm}[!t]
\begin{algorithmic}
\setstretch{1.2}
    \Function{loads}{$\bm{g},M_{start},t_f,\varepsilon_\nu,N_{max,e}\forall e$}
        \State Initialize $W_{t}=M_{start}$ and $M_f=\emptyset$.
        \While{$W_{t}\neq\emptyset$}
            \State Remove the first set $[\bm{x}^{(k)}(t_k)]$ from $W_{t}$.
            \While{$t_k<t_f$}
                \State Perform one integration step to obtain $[\bm{x}^{(k)}(t_k+\Delta t)]$.
                \State Compute the \acrshort{stm} $\left[\bm{\Phi}^{(k)}(t_k+\Delta t,t_0)\right]$.\Comment{see \cref{eq:STM}}
                \State Extract the non-constant part of the \acrshort{stm} $\delta\bm{\Phi}^{(k)}(t_k+\Delta t,t_0)$.\Comment{see \cref{eq:dSTM}}
                \State Compute the matrix of upper polynomial bounds $\bm{B}_{\delta\bm{\Phi}^{(k)}(t_k+\Delta t,t_0)}$.\Comment{see \cref{eq:dSTM_bounds}}
                \State Compute the \acrshort{nli} of the transformation $\nu^{(k)}(t_k+\Delta t,t_0)$.\Comment{see \cref{eq:nu_STM}}
                \If{$\nu^{(k)}(t_k+\Delta t,t_0)\leq\varepsilon_\nu$}
                    \If{$t_k+\Delta t=t_f$}
                        \State Store $[\bm{x}^{(k)}(t_k+\Delta t)]$ into $M_f$.
                    \Else
                        \State Update $t_k=t_k+\Delta t$.
                    \EndIf
                \Else
                    \State Compute the splitting direction $d$.\Comment{see \cref{eq:nu_e_STM,eq:dir_STM}}
                    \State Compute the number of splits $N_d^{(k)}$ along $d$.
                    \If{$N_d^{(k)}<N_{max,d}$}
                        \State Split $[\bm{x}^{(k)}(t_k)]$ into the $[\bm{x}^{(k+s)}(t_k)]$ sets.\Comment{see \cref{eq:LOADS_compose}}
                        \State Store $\{[\bm{x}^{(k+s)}(t_k)]\}$ into $W_{t}$.
                    \Else
                        \State Propagate $[\bm{x}^{(k)}(t_k)]$ to $t_f$ without splits.
                        \State Store $[\bm{x}^{(k)}(t_f)]$ into $M_f$.
                    \EndIf
                \EndIf
            \EndWhile
        \EndWhile
        \State \Return $M_{f}$
    \EndFunction
\end{algorithmic}
\caption{\glsentryshort{loads} algorithm for uncertainty propagation}
\label{alg:loads_algorithm_up}
\end{algorithm}
\hl{The \acrshort{nli} of the transformation is finally given by}
\begin{equation}
\label{eq:nu_STM}
\begin{aligned}
\nu^{(k)}(t_k+\Delta t,t_0) &= \sqrt{\dfrac{\displaystyle\sum_{i=1}^m\sum\limits_{j=1}^n\left(\sum\limits_{p=1}^n\left\lvert c_{ij,p}^{(k)}\right\rvert\right)^2}{\displaystyle\sum\limits_{i=1}^m\sum\limits_{j=1}^n \left(\bar{\Phi}^{(k)}_{ij}(t_k+\Delta t,t_0)\right)^2}}
\end{aligned}
\end{equation}
If the value of $\nu$ is lower than the imposed tolerance, the propagation is resumed, and continues until the nonlinearity constraints are satisfied or the final propagation epoch $t_f$ is reached. Conversely, if $\nu\hl{>}\varepsilon_{\nu}$, a splitting direction is computed by relying on the directional Jacobians, i.e.
\begin{equation}
\label{eq:nu_e_STM}
\nu_e^{(k)}(t_k+\Delta t,t_0) = \sqrt{\dfrac{\displaystyle\sum_{i=1}^m\sum\limits_{j=1}^n\left(c^{(k)}_{ij,e}\right)^2}{\displaystyle\sum\limits_{i=1}^m\sum\limits_{j=1}^n \left(\bar{\Phi}^{(k)}_{ij}(t_k+\Delta t,t_0)\right)^2}}\qquad e=1,\ldots,n
\end{equation}
\begin{equation}
\label{eq:dir_STM}
d = \argmax_{e}\left\{\nu^{(k)}_{e}(t_k+\Delta t,t_0)\right\}
\end{equation}
At this point, if the number of splits along $d$ is lower than the maximum allowed, the set $\left[\bm{x}^{(k)}(t_k)\right]$ is split into three new sets, thus obtaining
\hl{
\begin{subequations}
\begin{gather}
    \delta\bm{x}_{d}^{(k+s)} = \left\{\delta x_1,\ldots,\delta x_{d-1},\dfrac{1}{3}\delta x_d+\dfrac{2}{3}(s-1),\delta x_{d+1},\ldots,\delta x_n\right\}^{\textrm{T}},\qquad s=\{0,1,2\}\\
    \left[\bm{x}^{(k+s)}(t_k)\right] = \left[\bm{x}^{(k)}(t_k)\right]\circ\delta\bm{x}^{(k+s)}_{d}\\
    \lambda^{(k+s)}_{d} = \dfrac{1}{9}\lambda^{(k)}_d
\end{gather}
\label{eq:LOADS_compose}
\end{subequations}
}
These sets replace set $k$ in the manifold $W_t$. Conversely, in case a split in a forbidden direction is required, \hl{i.e. $N_d^{(k)}=N_{max,d}$}, the set $\left[\bm{x}^{(k)}(t_k)\right]$ is propagated to the final epoch with no more splits allowed. The procedure \hl{is then repeated with} the first set of $W_t$, and terminates \hl{when} all the domains are processed\hl{, i.e. $W_t=\emptyset$. The output is the uncertainty set propagated at $t_f$ described by a manifold $M_f$ of second-order expansions generated by the \acrshort{loads} algorithm. Note that, even though a single split is performed each time the \acrshort{nli} threshold is violated, multiple splits can occur at the same propagation time $t_k$. When nonlinearities are detected in $[\bm{x}^{(k)}(t_k+\Delta t)]$, the last integration step is in fact rejected and the split is performed on $[\bm{x}^{(k)}(t_k)]$ to generate the $[\bm{x}^{(k+s)}(t_k)]$ that will populate $W_t$. When one of the $[\bm{x}^{(k+s)}(t_k)]$ is extracted from $W_t$, the polynomial state is propagated to $t_k+\Delta t$\footnote{\hl{$\Delta t$ is the integration step size chosen by the adaptive step size integrator. It is generally different for each step and processed set.}} and a new check is performed on the transformed state. If $\nu$ is still violated, the step is rejected and a second split is triggered again at $t_k$.}

\subsection{Merging breaks}
\label{subsec:merge}
\begin{figure}[]
\centering
\subfloat[\label{subfig:Merging_1}]{\includegraphics[trim=1cm 4.5cm 0cm 2cm, clip=true, width=0.72\textwidth]{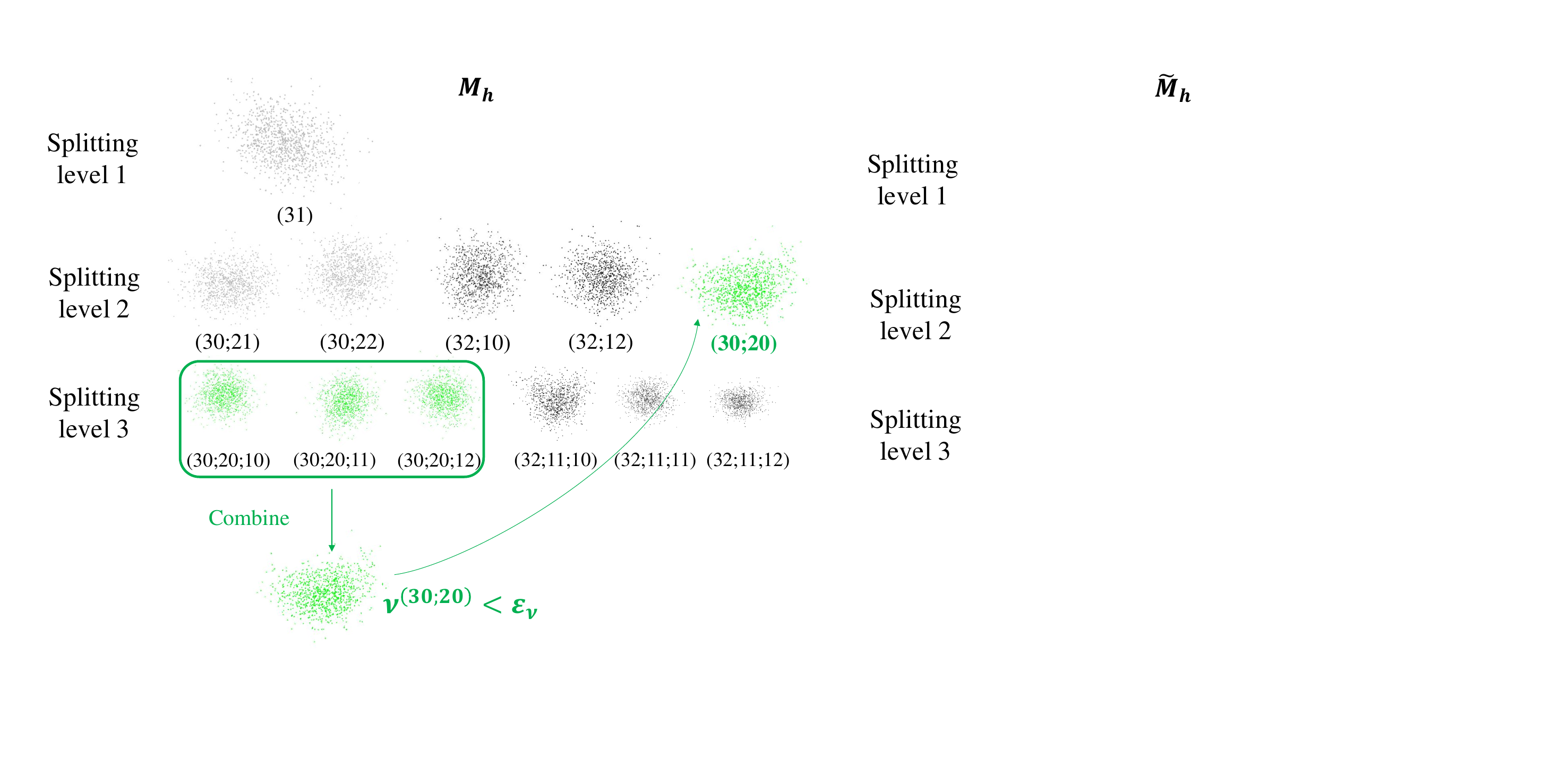}}\\
\subfloat[\label{subfig:Merging_2}]{\includegraphics[trim=1cm 4.5cm 0cm 2cm, clip=true, width=0.72\textwidth]{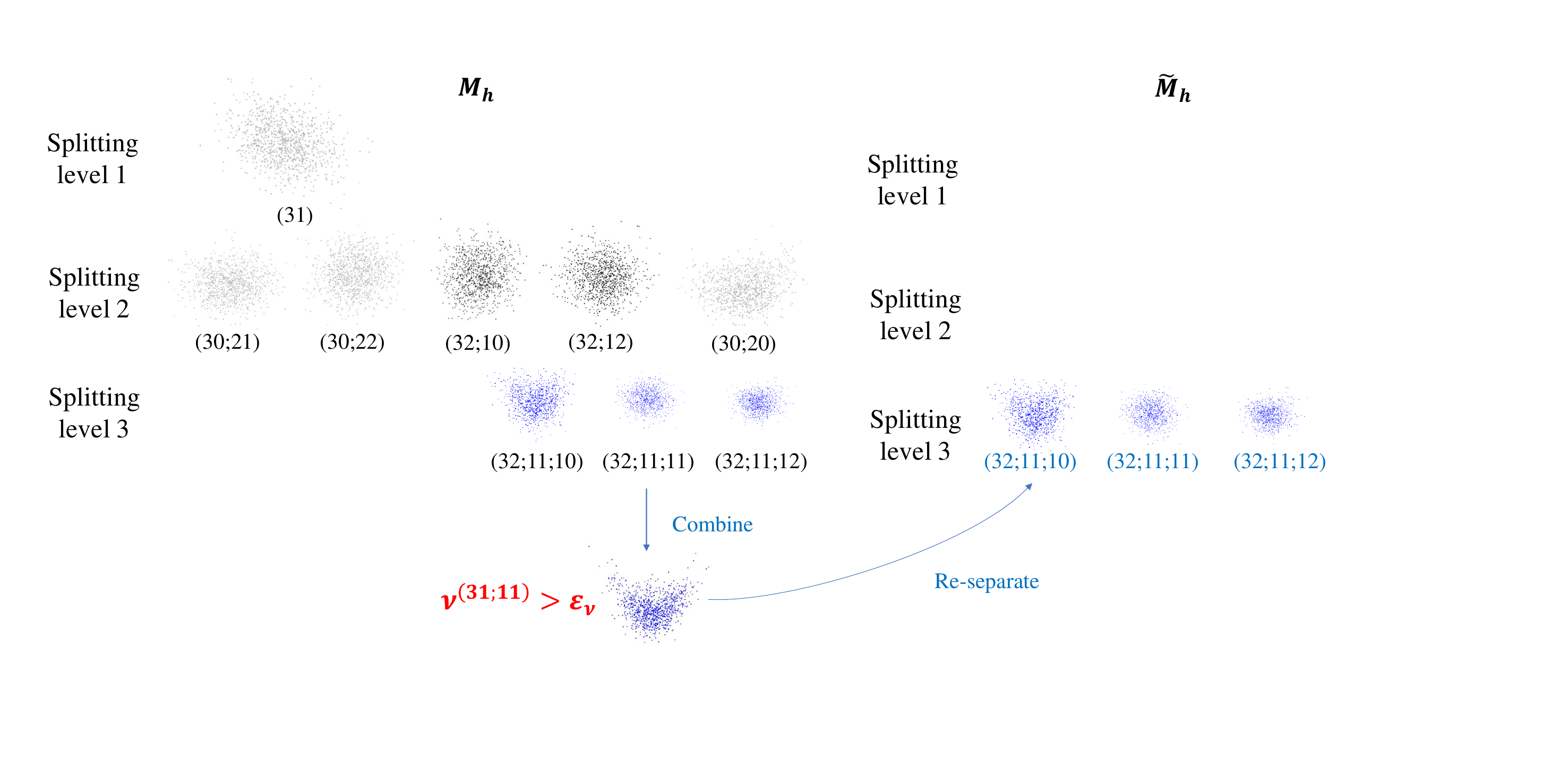}}\\
\subfloat[\label{subfig:Merging_3}]{\includegraphics[trim=1cm 4.5cm 0cm 2cm, clip=true, width=0.72\textwidth]{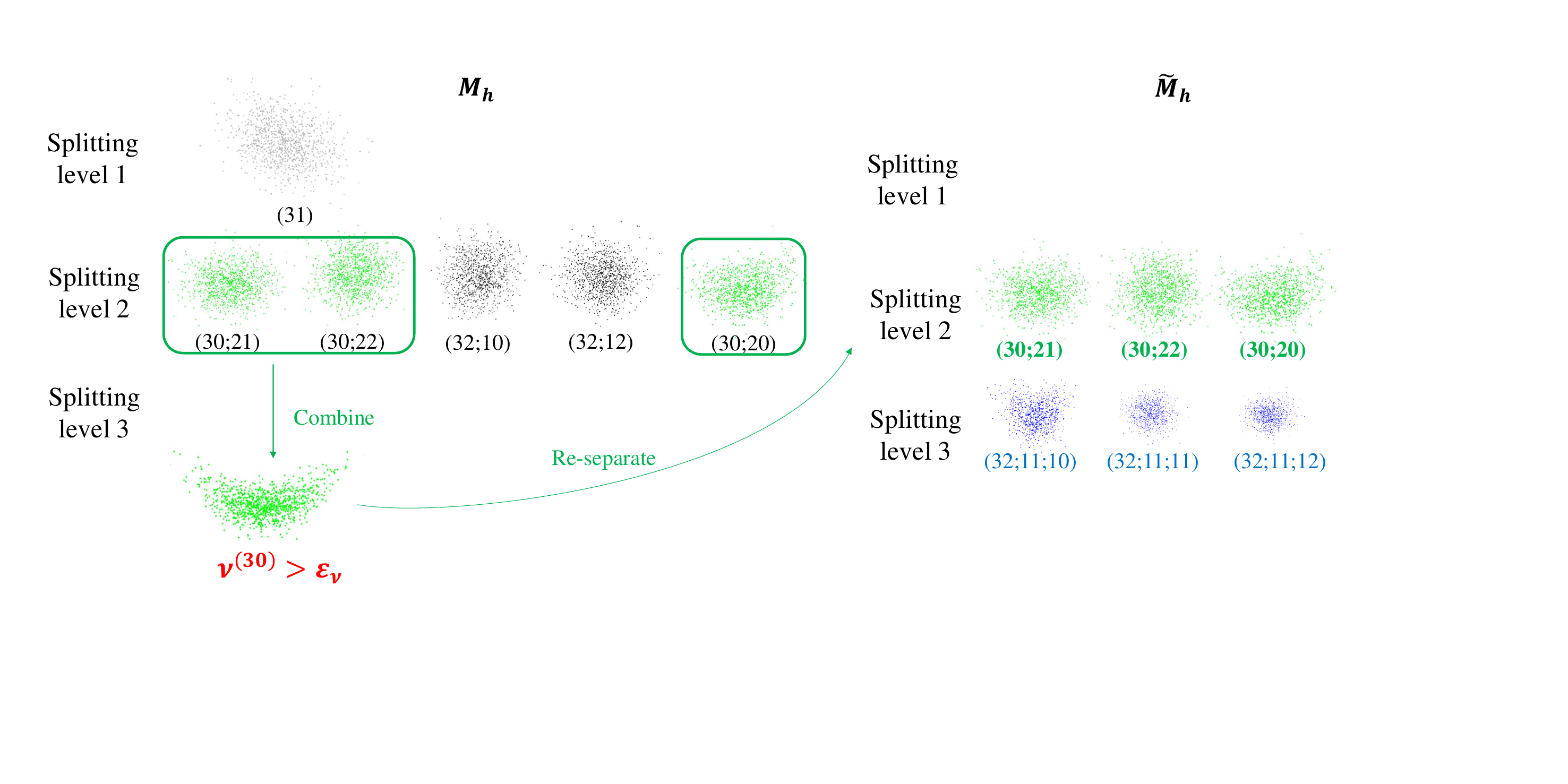}}\\
\subfloat[\label{subfig:Merging_4}]{\includegraphics[trim=1cm 4.5cm 0cm 2cm, clip=true, width=0.72\textwidth]{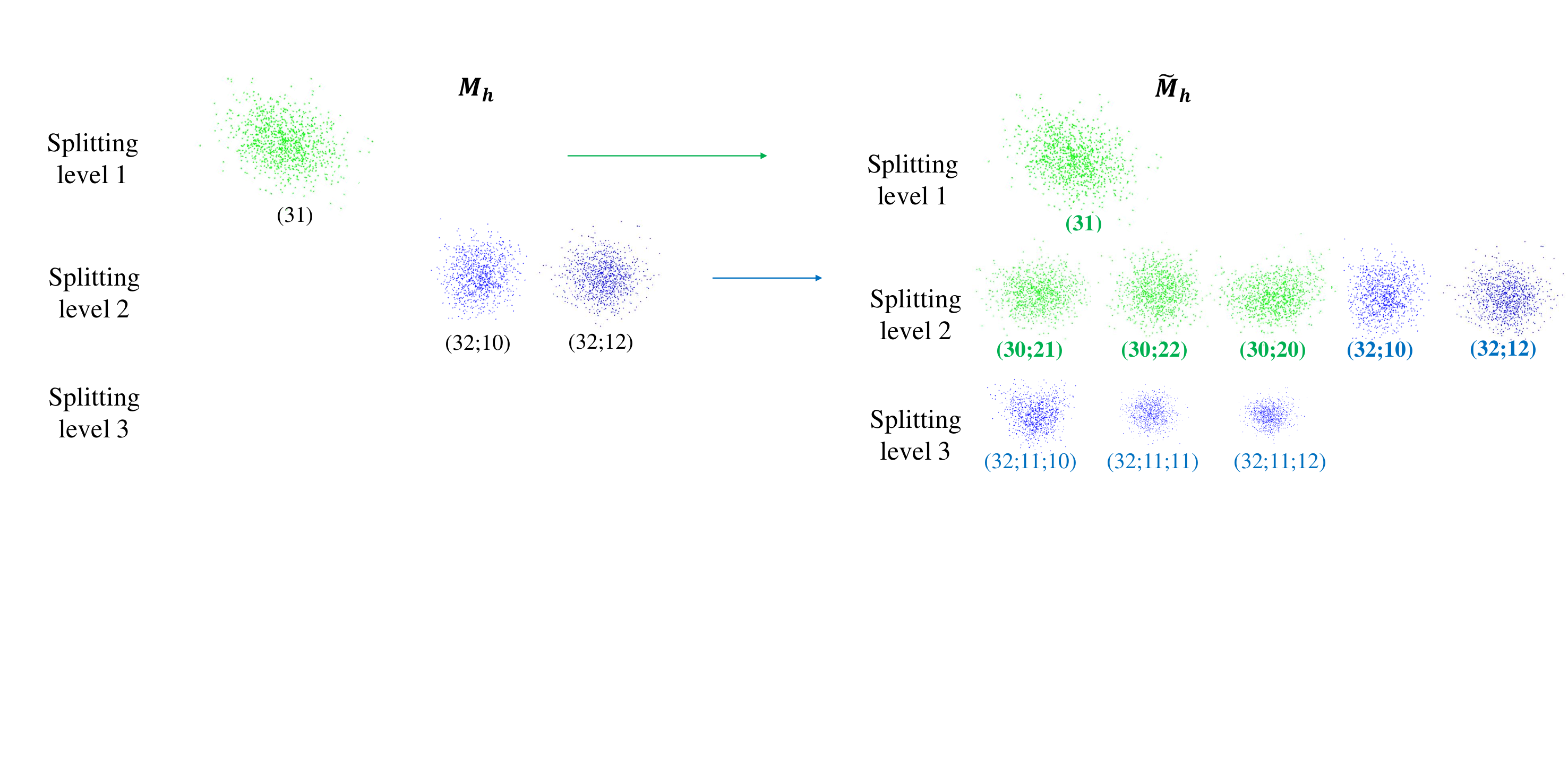}}\\
\caption{Merging phase description: \protect\subref{subfig:Merging_1} recombination of three sets at splitting level 3; \protect\subref{subfig:Merging_2} merging failure at level 3, and $\tilde{M}_h$ update; \protect\subref{subfig:Merging_3} merging failure at level 2 and $\tilde{M}_h$ update; \protect\subref{subfig:Merging_4} final $\tilde{M}_h$ composition.}
\label{fig:Merging}
\end{figure}

The number of sets generated by the \acrshort{loads} algorithm is governed by the nonlinearity levels encountered by each set during its propagation. The described approach is nonreversible, i.e. once generated, the sets are maintained throughout the propagation until further sets are required. This is actually a limitation since the nonlinearities can reduce along with the propagation. More specifically, nonlinearities are stronger at the pericenter and weaker at the apocenter. Consequently, the number of sets grows close to the pericenter while no extra split\hl{s are} expected far from it. Starting from these considerations, the proposed \acrshort{loads} algorithm for orbital uncertainty propagation is enhanced\hl{, resulting in the} split/merge algorithm described hereafter. In order to help the visualisation of the results, we \hl{now use} the splitting history notation introduced in Section~\ref{sec:LOADS}. As a result, a generic set $[\bm{x}^{(k)}]$ can be also written as  $[\bm{x}^{\left(H^{\left(k\right)}\right)}]$, i.e. can be univocally characterised by its splitting history. For example, a set ``(30;22;11)'' indicates a domain obtained from $[\bm{x}_0]$ by first splitting along direction 3 and choosing the left subset (``(30)''), then splitting  ``(30)'' along 2 and selecting the right one (``(30;22)''), and finally dividing ``(30;22)'' along direction 1 and selecting the middle set.

The merging algorithm makes use of some intermediate phases defined ``merging breaks'' that are introduced during the propagation. When a merging break is introduced at an epoch $t_h<t_f$, an intermediate manifold $M_{h}$ is built by propagating all sets to $t_h$, i.e.
\begin{equation}
M_h=\left\{\left[\bm{x}^{(k)}(t_k)\right]: t_k= t_h\right\}\qquad k=1\ldots,n_{set}
\end{equation}
The \hl{merging} algorithm is then run to reduce the size of $M_{h}$, thus recombining into a reduced manifold $\tilde{M}_h$ sets generated during the pericenter passage that are no longer needed. A summary of the procedure is given in Algorithm~\ref{alg:merge_algorithm}.

\begin{algorithm}[!t]
\begin{algorithmic}
\setstretch{1.2}
    \Function{\hl{Merging}}{$M_h,\varepsilon_\nu$}
        \State Group all sets in $M_h$ by depth of split to obtain $\mathcal{M}_h=\left\{M_h^0,\ldots,M_h^{N_{max}}\right\}$.
        \State Initialize $\tilde{M}_h=\emptyset$.
        \For{$z={N_{max}},{N_{max}}-1,\ldots,1$}
            \State Remove $M_h^z$ from $\mathcal{M}_h$.
            \While{$M_h^z\neq\emptyset$}
                \State Remove the first set from $M_h^z$.
                \State Find the siblings sets (if any) and remove them from $M_h^z$.
                \State Build the list $T$ with the found sets.
                \If{$\abs{T}=3$}
                    \State Identify the middle set of the triplet
                    \State Re-expand it.\Comment{See \cref{eq:merge_compose}}
                    \State Compute its \glsentryshort{nli}.\Comment{See \cref{eq:stm_dev_compose,eq:nli_compose}}
                    \If{$\nu\leq\varepsilon_\nu$}
                        \State Append the expanded set to $M_h^{z-1}$.
                    \Else
                        \State  Append the elements in $T$ to $\tilde{M}_h$.
                    \EndIf
                \Else
                    \State Append the elements in $T$ to $\tilde{M}_h$.
                \EndIf
            \EndWhile
        \EndFor
        \If{$M_h^0\neq\emptyset$}
            \State Set $\tilde{M}_h=M_h^0$.
        \EndIf
        \State \Return $\tilde{M}_h$.
    \EndFunction
\end{algorithmic}
\caption{\hl{Merging} algorithm}
\label{alg:merge_algorithm}
\end{algorithm}

The merging process is initialized by grouping the sets of $M_{h}$ into so-called ``splitting levels'' by considering their depth of split, i.e. how many splits they were subject to, thus building a list of manifolds $\mathcal{M}_h=\left\{M_h^0,M_h^1,\ldots,M_h^{N_{max}}\right\}$ with ${N_{max}}$ maximum depth of split among sets in $M_h$. The notation $M_h^{z}$ indicates a subset of manifold $M_h$ collecting all the domains that went through \hl{exactly} $z$ splits. Therefore  $M_h^{1}$ is the splitting level 1, i.e. the group of sets obtained from the initial domain with just one split.  Figure~\ref{fig:Merging} shows an example of how the sets can be grouped by depth of split into splitting levels. For instance, sets ``(30;20;10)'', ``(32;20;11)'', ``(30;20;12)'' and ``(32;11;12)'' are all characterised by three pairs, which means that they were obtained from the original set with three splits in sequence, thus they belong to level 3. Among sets with the same depth of split $N$, so called ``siblings sets'' can be recognized. These sets are triplets of domains that share the same splitting direction for all the couples of their splitting history, and the same splitting side for the first $(N-1)$ couples.
In our example sets ``(30;20;10)'', ``(30;20;11)'' and ``(30;20;12)'' are siblings sets. This is not true for set ``(32;11;12)'', which was obtained with a 3-1-1 splitting history. If a triplet of ``siblings sets'' is identified, a recombination is attempted. This is performed by re-expanding the middle set around its centre point. For example, consider the sets $\left[\bm{x}^{(30;20;10)}(t_h)\right]$, $\left[\bm{x}^{(30;20;11)}(t_h)\right]$ and $\left[\bm{x}^{(30;20;12)}(t_h)\right]$ of \cref{subfig:Merging_1}. The three sets are combined, thus re-obtaining the set $\left[\bm{x}^{(30;20)}(t_h)\right]$. The operation is performed by function composition, i.e., once defined the vector \hl{$\delta\bm{x}_1=\left\{3\delta x_1,\delta x_2,\ldots,\delta x_n\right\}^{\textrm{T}}$}, the new set can be obtained as
\begin{equation}
\label{eq:merge_compose}
\left[\bm{x}^{(30;20)}(t_h)\right] = \left[\bm{x}^{(30;20;11)}(t_h)\right]\circ\delta\bm{x}_1\qquad \lambda^{(30;20)}_1 = 9\lambda^{(30;20;11)}_1
\end{equation}
At this point, the new \acrshort{nli} can be computed as
\begin{equation}
\delta\bm{\Phi}^{(30;20)}(t_h,t_0) = \left[\dfrac{\partial\left[x^{(30;20)}_i\right]}{\alpha_{\lambda_j}\sqrt{\lambda^{(30;20)}_j}\partial\delta x_j}(\delta\bm{x})\right] - \left[\dfrac{\partial\left[x^{(30;20)}_i\right]}{\alpha_{\lambda_j}\sqrt{\lambda^{(30;20)}_j}\partial\delta x_j}(\bm{0})\right]= \left[\displaystyle\sum\limits_{p=1}^n c_{ij,p}^{(30;20)}\delta x_p\right]
\label{eq:stm_dev_compose}
\end{equation}
\begin{equation}
\nu^{(30;20)}(t_h,t_0) =\sqrt{\dfrac{\displaystyle\sum\limits_{i=1}^m\sum\limits_{j=1}^n\left(\sum\limits_{p=1}^n\left\lvert c_{ij,p}^{(30;20)}\right\rvert\right)^2}{\displaystyle\sum\limits_{i=1}^m\sum\limits_{j=1}^n\left(\bar{\Phi}_{ij}^{(30;20)}(t_h,t_0)\right)^2}}
\label{eq:nli_compose}
\end{equation}
and compared with the imposed threshold $\varepsilon_{\nu}$. If the accuracy requirement is satisfied, the recombined set is retained and added to the splitting level 2. Otherwise, the recombined set is discarded and the sibling sets appended to a reduced manifold $\tilde{M}_h$. The procedure continues for all possible triplets of the current level. Then, the lower level is investigated, and the whole process is repeated until $\mathcal{\bm{M}}_h$ is emptied. 

An example of the described process is shown in Fig.~\ref{fig:Merging}. Three splitting levels are identified, meaning that all sets were split three times at maximum. At splitting level 3, sets ``(30;20;10)'',``(30;20;11)'' and ``(30;20;12)'' are successfully recombined, and the resulting merged set is saved as ``(30;20)'' at splitting level 2. The triplet ``(32;11;10)'',``(32;11;11)'' and ``(32;11;12)'' instead does not satisfy the linearity threshold when recombined. These sets are thus directly added to $\tilde{M}_h$. At this point, the analysis passes to the splitting level 2. Sets ``(30;21)'' and ``(30;22)'' are combined with the just merged ``(30;20)'' set. Their recombination, however, does not respect the imposed threshold $\varepsilon_{\nu}$; thus, the three sets are added to level 2 of the reduced manifold. At this point, no further recombination is possible, and the remaining three sets (sets ``(32;10)'' and ``(32;12)'' at level 2 and set ``(31)'' at level 1) are automatically moved to $\tilde{M}_h$. The recursive process is summarized in \cref{alg:merge_algorithm}. The advantage of the merging break is evident: the initial manifold was made of 11 sets, which reduce to 9 after the merging. Once the merging phase is over, the propagation is resumed from $\tilde{M}_h$.


\subsection{Combined \acrshort{loads}-\hl{merging} algorithm}
\begin{algorithm}[!t]
\begin{algorithmic}
\setstretch{1.2}
    \Function{\acrshort{loads}\hl{Merging}}{$\bm{g},[\bm{x}(t_0)],t_0,T,n_r,n_{b/r},\varepsilon_\nu,N_{max,e}\forall e$}
        \State Build the epochs vector $\bm{t}$.\Comment{See \cref{eq:dt_merge}}
        \State Initialize $\tilde{M}_{0}=\{[\bm{x}(t_0)]\}$.
        \For{$h=1,\ldots,\abs{\bm{t}}$}
            \State Run \cref{alg:loads_algorithm_up} with $M_{start}=\tilde{M}_{h-1}$ and $t_f=t_h$ to obtain $M_{h}$.
            \State Run \cref{alg:merge_algorithm} with $M_{h}$ to obtain $\tilde{M}_h$.
        \EndFor
        \State \Return $\tilde{M}_{\abs{\bm{t}}}$.
    \EndFunction
\end{algorithmic}
\caption{\acrshort{loads}-\hl{merging} algorithm for orbital uncertainty propagation.}
\label{alg:loads_merge_algorithm}
\end{algorithm}
The algorithm in \cref{subsec:merge} is combined with the splitting scheme as follows. A number of merging breaks per revolution $n_{b/r}$ is firstly specified and an epoch vector $\bm{t}=\{t_0,\ldots,t_q\}$ is built such that
\begin{equation}
    t_h = t_{h-1} + \dfrac{n_r\cdot T}{\lceil n_{b/r}\cdot n_r \rceil}\qquad \forall\ h\in[1,q]
    \label{eq:dt_merge}
\end{equation}
with $t_0$ initial time, $T$ orbit period, $n_r$ number of revolutions and $q=\lceil n_{b/r}\cdot n_r \rceil$. Note that $n_r$ might be a fractional number. The \acrshort{loads} algorithm is then run for each time interval $[t_{h-1},t_{h}]$ followed by a \hl{merging} attempt at $t_h\ \forall h\in[1,q]$ thus resulting in the recursive scheme summarized in \cref{alg:loads_merge_algorithm}.

\section{LOADS-GMM}
\label{sec:LOADS-GMM}
The described \acrshort{loads} algorithm provides a manifold of polynomials \hl{$M_{\bm{y}}$ \hl{(or $M_{f}$, in case of orbital uncertainty propagation)} that maps $[\bm{x}]$} through the considered nonlinear transformation. \hl{Recalling that $[\bm{x}]$ is the \acrshort{da} representation of the random variable $\bm{X}$ with \acrshort{pdf} $p_{\bm{X}}$, this} map can be proficiently used to obtain a statistical characterisation of the propagated \hl{random} variable \hl{$\bm{Y}$}.
More specifically, \hl{realizations $\bm{x}_l$ of $\bm{X}$ can be drawn for $l=1,\ldots,N_{MC}$ by standard \acrshort{mc} sampling, where $N_{MC}$ is the total number of samples.} Each realization is then associated to a single domain \hl{in} $M_{\bm{y}}$ (or $M_{f}$), and then mapped \hl{to $\hat{\bm{y}}_l$} via polynomial evaluation. Once all realizations are mapped, higher order moments of the \hl{transformed \acrshort{pdf} $p_{\bm{Y}}$} can be computed \hl{from $\hat{\bm{y}}_l$}~\citep{Pebay2016}. This approach allows us to statistically describe the propagated uncertainty, but does not provide an analytical formulation of its \acrshort{pdf}. The availability of such an analytical representation can be appealing in a variety of applications, including state estimation \cite{Horwood2011}, tracks correlations \cite{DeMars2013a}, collision probability computation \cite{DeMars2014a}, and collision avoidance design \cite{Armellin2021347}.

While the problem of estimating the \acrshort{pdf} of a Gaussian distribution subject to a linear transformation is of trivial solution since the Gaussian nature is preserved, the mapping of a generic \acrshort{pdf} through a nonlinear transformation \hl{$\bm{f}$} is more complex. A possible approach relies on the use of so-called \acrfullpl{gmm}. According to the \acrshort{gmm} theory, any \acrshort{pdf} \hl{of practical concern} can be approximated by a properly selected weighted sum of Gaussian distributions or kernels~\citep{Sorenson1971}. As a result, given a distribution $p_{\bm{X}}$, we can write
\begin{equation}
\label{eq:p_x}
p_{\bm{X}}\approx \sum_{q=1}^{N_x}w_{\bm{X}}^{(q)} p_g\left(\bm{x};\bm{\mu}_{\bm{X}}^{(q)},\bm{P}_{\bm{X}}^{(q)}\right)
\end{equation}
where $p_g$ indicates a generic Gaussian \acrshort{pdf}, $\bm{\mu}_{\bm{X}}^{(q)}$ and $\bm{P}_{\bm{X}}^{(q)}$ are the mean and the covariance matrix of the kernel $q$, whereas $w_{\bm{X}}^{(q)}\geq 0$ \hl{are} the associated weight\hl{s} such that $\sum_q w_{\bm{X}}^{(q)} =1$. Now suppose that $p_{\bm{X}}$ is known. This may coincide with a single Gaussian kernel, or may be the result of an estimation process. Consider the generic nonlinear transformation $\bm{y}=\bm{f}(\bm{x})$ and assume that we want to compute $p_{\bm{Y}}$. Without loss of generality, one can always write
\begin{equation}
\label{eq:p_y}
p_{\bm{Y}}\approx \sum_{p=1}^{N_y}w_{\bm{Y}}^{(p)} p_g\left(\bm{y};\bm{\mu}_{\bm{Y}}^{(p)},\bm{P}_{\bm{Y}}^{(p)}\right)
\end{equation}
A direct transformation from \cref{eq:p_x} to \cref{eq:p_y} is generally not straightforward. However, if we tune the number and shape of the components of $p_{\bm{X}}$ such that $\bm{f}$ is locally linear or quasi-linear on a proper domain around each kernel mean, several methods exist to propagate each Gaussian kernel of $p_{\bm{X}}$ through $\bm{f}$. Then, assuming the weights of the \acrshort{gmm} can be maintained constant~\citep{Horwood2011}, the propagated kernels can be seen as the components of a new \acrshort{gmm} approximating the random variable $\bm{Y}$. That is, \cref{eq:p_y} is automatically retrieved.

Such a \acrshort{gmm}-based uncertainty propagation approach requires two ingredients, i.e. a method for splitting Gaussian kernels, and a criterion for deciding how to split them. A well-established literature on splitting routines exists, and it is based on so-called splitting libraries. Following the description given by~\citet{DeMars2013}, let us first start from a 1D case, and consider a standard Gaussian distribution $p(x)=p_g(x;0,1)$. The goal of the GMM method is to find a $\tilde{p}(x)$ such that
\begin{equation}
p(x)\approx\tilde{p}(x) = \sum\limits_{s=0}^{L-1}\tilde{w}^{(s)} p_g(x;\tilde{\mu}^{(s)},\tilde{\sigma}^{(s)})
\end{equation}
where $L$ is the number of components, $\tilde{\mu}^{(s)}$ and $\tilde{\sigma}^{(s)}$ are the mean and standard deviation of the $s$-th Gaussian kernel, respectively, whereas $\tilde{w}^{(s)}\geq 0$ \hl{are} the associated weight\hl{s} such that $\sum_s\tilde{w}^{(s)}=1$. The solution is found with an optimization problem aiming at minimizing a properly defined distance between $p$ and $\tilde{p}$. In the adopted approach, the $L_2$ distance is selected, whereas all $\tilde{\sigma}^{(s)}$ are set equal to a single value $\tilde{\sigma}$ and contribute to the cost function through a term $\lambda\tilde{\sigma}$, with $\lambda\in[0,1]$. As a result, so-called univariate splitting libraries are obtained. These libraries provide the values of $\tilde{w}^{(s)}$, $\tilde{\mu}^{(s)}$ and $\tilde{\sigma}$ as a function of the number of components $L$ and the parameter $\lambda$. The latter governs the value of $\tilde{\sigma}$, i.e. $\lambda\to 1$ leads to smaller variances and $\tilde{\mu}^{(s)}$ farther from zero while the opposite holds for $\lambda\to 0$. Examples can be found in~\citet{DeMars2013} and~\citet{Vittaldev2016a}. In this work, splitting libraries are computed numerically by using the interior point optimizer IPOPT~\citep{Wachter2006} compiled against the HSL Mathematical Software Library for the solution of sparse linear systems~\citep{HSL2013}.

Let us now consider an $n$-dimensional case, and assume we want to replace a single component of the GMM, that is
\begin{equation}
\label{eq:GMM_n}
w p_g(\bm{x};\bm{\mu},\bm{P})\approx\sum\limits_{s=0}^{L-1}w^{(s)} p_g\left(\bm{x};\bm{\mu}^{(s)},\bm{P}^{(s)}\right)
\end{equation}
A possible solution for the problem is presented in~\citet{DeMars2013} and exploits the eigendecomposition of matrix $\bm{P}$. As a result, assuming a split along the $d$-th eigenvector, the unknowns of Eq.~\eqref{eq:GMM_n} can \hl{be} found as
\begin{equation}
\label{eq:jump}
w^{(s)} = w\tilde{w}^{(s)}\qquad\bm{\mu}^{(s)}=\bm{\mu}+\sqrt{\lambda_d}\tilde{\mu}^{(s)}\bm{v}_d\qquad\bm{P}^{(s)}=\bm{V}\tilde{\bm{\Lambda}}\bm{V}^{\textrm{T}}
\end{equation}
where $\bm{v}_d$ is the $d$-th eigenvector of $\bm{P}$, $\bm{V}$ the eigenvector matrix, and $\tilde{\bm{\Lambda}}=\textrm{diag}\left\{\lambda_1,\ldots,\lambda_{d-1},\tilde{\sigma}^2\lambda_d,\lambda_{d+1},\ldots,\lambda_n\right\}$, with $\{\lambda_1,\ldots,\lambda_n\}$ eigenvalues of $\bm{P}$.

The above expressions illustrate how to arbitrarily split a generic kernel of a given \acrshort{gmm} distribution. Still, a method to decide how to split the initial \acrshort{gmm} to later retrieve the propagated \acrshort{pdf} is missing. Two possible philosophies can be followed. One may decide to a priori generate a large number of kernels~\citep{Vittaldev2016a}, or exploit the known nonlinear transformation to adaptively trigger splits~\citep{DeMars2013}, adjust the components weights~\citep{Terejanu2008}, and merge components when required~\citep{Terejanu2011,Vishwajeet2018}. Adaptive approaches have the advantage of exploiting the knowledge of the nonlinear transformation, but typically set up an augmented problem, as additional equations are required to estimate the nonlinearity of the problem. As an example, \citet{DeMars2013} require the estimation of the linearized and nonlinear differential entropy, which imply the numerical integration of additional equations, while~\citet{Terejanu2008} minimize the integral square difference between the true forecast \acrshort{pdf} and its \acrshort{gmm} approximation.

The approach we propose belongs to the family of adaptive approaches and has the advantage of naturally exploiting the result of the \acrshort{loads} algorithm. As previously described, the method automatically divides the initial uncertainty set in order to grant that over each generated domain the linearity assumption locally holds.
If we force the splitting algorithm to apply univariate splitting libraries whenever splits are required, and then we attach a Gaussian kernel to each generated subset, then we are automatically building the \acrshort{gmm} approximation of $p_{\bm{X}}$ we need. By definition of \acrshort{nli}, each generated subset is such that the behaviour of $\bm{f}$ on it can be assumed as linear or quasi-linear. Polynomial evaluation on \acrshort{ut} sigma points is then used to propagate each Gaussian kernel through $\bm{f}$. The propagated \acrshortpl{pdf} will then constitute the kernels of the \acrshort{gmm} approximation of $p_{\bm{Y}}$. Unlike existing adaptive algorithms, the detection of nonlinearities and the identification of the splitting directions are automatically granted by the \acrshort{da}-based formulation of the problem, thus no solution of any further equation is required.

\begin{figure}[!t]
\centering
\includegraphics[trim=1.5cm 2cm 1.5cm 2cm, clip=true, width=\textwidth]{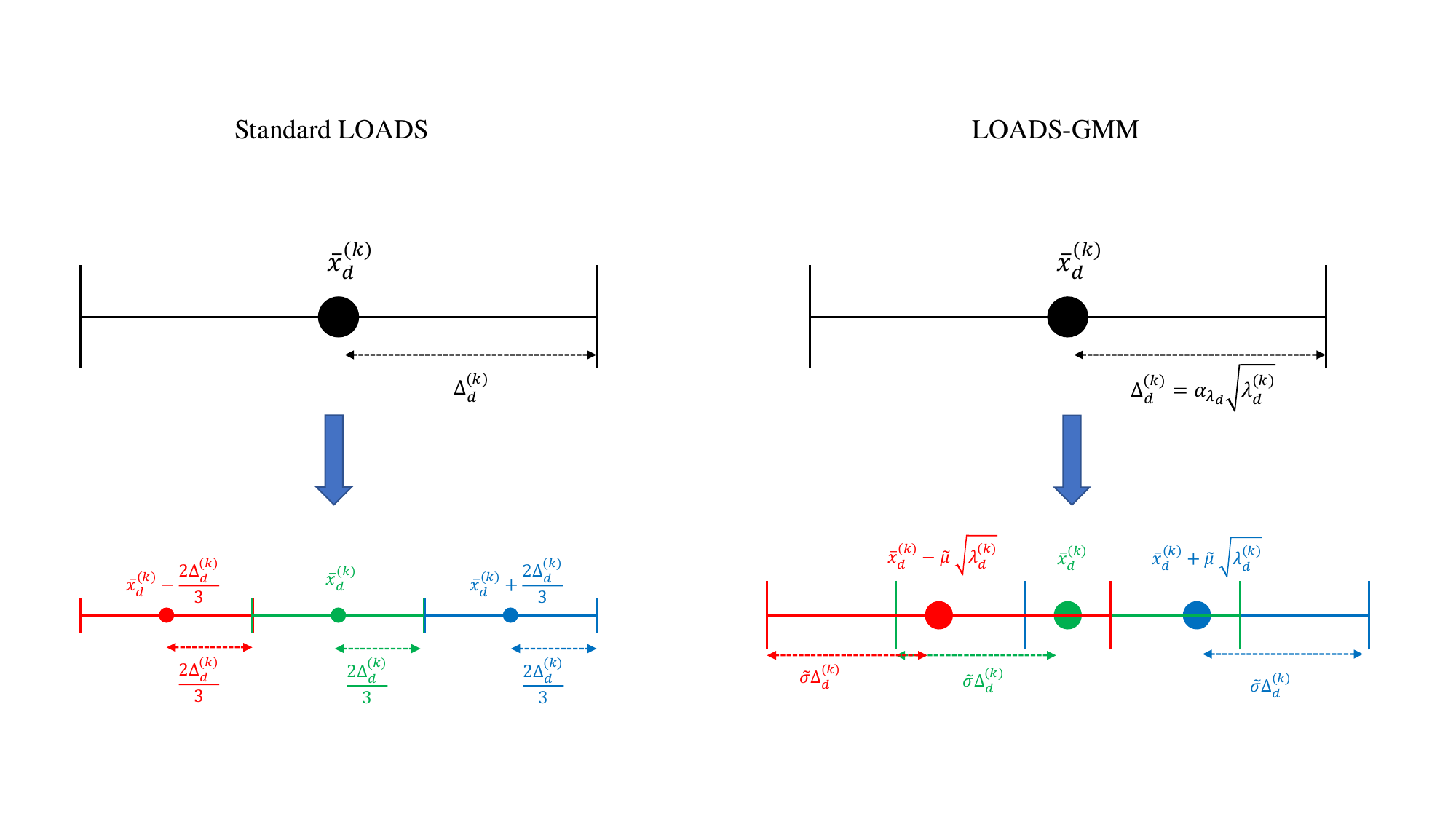}\\
\caption{Splitting routines for \acrshort{loads} and \acrshort{loads}-\acrshort{gmm}. On the left, the standard \acrshort{loads} subsets generation with 3 domains without overlapping. On the right, the splitting library used for \acrshort{gmm} generation ($L=3$, $\lambda = 1e-3$).}
\label{fig:SL}
\end{figure}

The resulting algorithm, here referred to as \acrshort{loads}-GMM, closely follows the standard \acrshort{loads} algorithm with some minor modifications, and it is here illustrated in the case of orbital uncertainty propagation. Under the assumptions of Section~\ref{sec:LOADS}, the process is initialized by setting $M_0=\left\{\left[\bm{x}_0\right]\right\}$, while a single component GMM is set up, namely
\begin{equation}
p_{\bm{X}_0}=p_g(\bm{x};\bm{\mu}_0,\bm{P}_{0})
\end{equation}
Then, the propagation starts. Let us now assume that, at a certain epoch $t$, a generic element $\left[\bm{x}^{(k)}(t)\right]$ of $W_t$ violates the nonlinearity threshold, thus requiring a split. The element $\left[\bm{x}^{(k)}(t)\right]$ contributes to $p_{\bm{X}_0}$ with its weight, mean and covariance, namely $w^{(k)}p_g\left(\bm{x};\bm{\mu}^{(k)},\bm{P}^{(k)}\right)$. The splitting procedure is similar to that of the standard \acrshort{loads}, but with a significant difference. The set, indeed, is no longer divided into 3 non-overlapping sets. Conversely, the univariate splitting library selected for building the GMM is used. A comparison between the two different approaches is illustrated in Fig.~\ref{fig:SL}. The plot shows the splitting along the direction $d$ of a generic $\left[\bm{x}^{(k)}\right]$. On the left, the standard \acrshort{loads} approach, with the generation of 3 subsets without overlapping. On the right, the application of a univariate splitting library. The selected library has $L=3$ and $\lambda=1e-3$. As a result, $\tilde{\sigma} = 0.6715664864669252$ whereas $\tilde{\mu}^{(1)} = 0$, while $\tilde{\mu}^{(0)} = \tilde{\mu}^{(2)}=\tilde{\mu}=1.0575150485760967$. As can be seen, the generated sets are bigger and show a significant overlap.

The generation of the subsets is performed by function composition. For example, assuming a split along direction $d$ and a splitting library $\left(\tilde{w}^{(s)},\tilde{\mu}^{(s)},\tilde{\sigma}\right)$ with $s=0,\ldots,L-1$, the function composition of Eq.~\eqref{eq:LOADS_compose} is replaced \hl{by}
\hl{
\begin{subequations}
\begin{gather}
    \delta\bm{x}_{d}^{(k+s)}=\left\{\delta x_1,\ldots,\delta x_{d-1},\tilde{\sigma}\delta x_d+\dfrac{\tilde{\mu}^{(s)}}{\alpha_{\lambda_d}},\delta x_{d+1},\ldots,\delta x_n\right\}^{\textrm{T}}\\
    \left[\bm{x}^{(k+s)}(t_k)\right] = \left[\bm{x}^{(k)}(t_k)\right]\circ \delta\bm{x}_{d}^{(k+s)}\\
    \lambda_{d}^{(k+s)} = \tilde{\sigma}^2\lambda_d^{(k)}
\end{gather}
\end{subequations}
}
The new generated sets $\left[\bm{x}^{(k+s)}(t)\right]$ are then coupled with weights, means, and covariances that are retrieved from those of $\left[\bm{x}^{(k)}(t)\right]$ as illustrated by Eq.~\eqref{eq:jump}, i.e.
\hl{
\begin{subequations}
\begin{gather}
    w^{(k+s)} = \tilde{w}^{(s)}w^{(k)}\\
    \bm{\mu}_0^{(k+s)} = \bm{\mu}_0^{(k)}+\sqrt{\lambda_d^{(k)}}\tilde{\mu}^{(s)}\bm{v}_d\\
    \bm{P}_0^{(k+s)} = \bm{V}\tilde{\bm{\Lambda}}^{(k)}\bm{V}^{\textrm{T}}
\end{gather}
\end{subequations}
}
Note that $\bm{v}$ and $\bm{V}$ have no apex since they coincide with the ones of $\bm{P}_0$. Splits, indeed, are performed in the set of the original coordinates, so no change of the eigenspace occurs. Once the sets are generated, they are added to the working list $W_t$, and the propagation is resumed.

The modification of the splitting library also affects the merging phase. The procedure is exactly the same as described in Section~\ref{subsec:merge}, but at each splitting level the algorithm now searches for $L$-tuples of sibling sets. As a result, if we indicate with $m=(L-1)/2$ the middle element of an $L$-tuple, a generic re-composition equation such as Eq.~\eqref{eq:merge_compose} can be re-expressed as
\begin{equation}
\left[\bm{x}^{(30;20)}(t_h)\right]=\left[\bm{x}^{(30;20;1m)}(t_h)\right]\circ\hl{\left\{\dfrac{\delta x_1}{\tilde{\sigma}},\delta x_2,\ldots,\delta x_n\right\}^{\textrm{T}}}\qquad \lambda_1^{(30;20)} = \dfrac{\lambda_1^{(30;20;1m)}}{\tilde{\sigma}^2}
\end{equation}
In addition, the \acrshort{pdf} elements of the recombined set can be computed as
\hl{
\begin{subequations}
\begin{gather}
    w^{(30;20)}=\dfrac{w^{(30;20;1m)}}{\tilde{w}^{(m)}}\\
    \bm{\mu}^{(30;20)}_0=\bm{\mu}_0^{(30;20;1m)}-\sqrt{\lambda_1^{(30;20)}}\tilde{\mu}^{(\hl{m})}\bm{v}_1\\
    \bm{P}^{(30;20)}_0 = \bm{V}\textrm{diag}\left\{\lambda_1^{(30;20)},\lambda_2^{(30;20;1m)},\ldots,\lambda_n^{(30;20;1m)}\right\}\bm{V}^{\textrm{T}}
\end{gather}
\end{subequations}
}
Once recombined the set, the procedure for the \acrshort{nli} computation is the one described in Section~\ref{subsec:merge}.

At the end of the process, the \acrshort{loads}-GMM algorithm provides two manifolds $M_0$ and $M_f$ of $n_{set}$ sets, and a GMM representation of the initial \acrshort{pdf} $p_{\bm{X}_0}$ which can be written as
\begin{equation}
p_{\bm{X}_0}\approx\sum\limits_{k=1}^{n_{set}} w^{(k)}p_g(\bm{x};\bm{\mu}_0^{(k)},\bm{P}_0^{(k)})
\end{equation}
The transformed \acrshort{pdf} at $t_f$ is then retrieved by individually propagating each GMM kernel. In the adopted approach, the \acrshort{ut} method is employed. Given a \acrshort{pdf}, the method selects \hl{a number of} so-called sigma points to accurately describe it. Each sigma point is then propagated independently, and a proper weighted sum of their images is used to reconstruct mean and covariance of the propagated \acrshort{pdf}. A detailed description of the method is given in~\citet{Julier2000}. In our formulation, given a generic element $k$ of $M_0$ and the corresponding \acrshort{pdf} parameters $\left\{w^{(k)},\bm{\mu}_0^{(k)},\bm{P}_0^{(k)}\right\}$, $2n$ samples $\bm{x}^{(k)}_{0,l}$, $l=1,\ldots,2n$ are drawn, with $n$ state space dimension. Each sample is then mapped to $t_f$ by exploiting the polynomial $\left[\bm{x}^{(k)}(t_f)\right]$, thus obtaining $\hat{\bm{x}}^{(k)}_{f,l}$. The first two moments of the propagated \acrshort{pdf} are then computed as
\begin{equation}
\begin{gathered}
\bm{\mu}_{f}^{(k)}=\sum\limits_{l=1}^{2n} \kappa_l\hat{\bm{x}}^{(k)}_{f,l}\\
\bm{P}_{f}^{(k)}=\sum\limits_{l=1}^{2n}\kappa_l\left(\hat{\bm{x}}^{(k)}_{f,l}-\bm{\mu}^{(k)}_{f}\right)\left(\hat{\bm{x}}^{(k)}_{f,l}-\bm{\mu}^{(k)}_{f}\right)^{\textrm{T}}
\end{gathered}
\end{equation}
where $\kappa_l$ is the \acrshort{ut} weight associated to sample $l$~\cite{Julier2000}. The procedure is then repeated for all the elements of $M_0$. At the end, the \acrshort{gmm} associated to the propagated \acrshort{pdf} $p_{\bm{X}_f}$ can be obtained as
\begin{equation}
p_{\bm{X}_f}\approx\sum\limits_{k=1}^{n_{set}}w^{(k)}p_g\left(\bm{x};\bm{\mu}^{(k)}_{f},\bm{P}^{(k)}_{f}\right)
\end{equation}

\section{Numerical simulations}
\label{sec:Sim}
This section illustrates some numerical applications of the mathematical concepts introduced in~\cref{sec:Index,sec:LOADS,sec:LOADS_unc,sec:LOADS-GMM}. The proposed \acrshort{da}-based \acrshort{nli} is compared with a reference \acrshort{nli} in Section~\ref{subsec:NLI_test}. The propagation of a large uncertainty set on a highly eccentric orbit is then used as a stress test to assess the performance of the \acrshort{loads} algorithm. Finally, the combined \acrshort{loads}-GMM method for orbital uncertainty propagation is directly applied on test cases available in literature. \hl{The presented methods are implemented in JAVA and exploit CNES' PACE library for Taylor algebra}. All the simulations presented in this paper were run on an Intel i7-8565U CPU @ 1.80GHz and 16GB of RAM.

\subsection{The \acrlong{nli}}
\label{subsec:NLI_test}
The performance of the \acrshort{da}-based \acrshort{nli} described in Section~\ref{sec:Index} is now tested and compared against a benchmark.
The reference method is the \acrshort{nli} formulated by~\citet{Junkins2004}, which inspired our approach. Given a nonlinear transformation $\bm{y} = \bm{f}(\bm{x})$ and a set $\bm{x}$ with nominal value and uncertainty, the index, here referred to as $\nu_{JS}$, is computed by sampling the worst case uncertainty surface of the considered problem, mapping each sample independently via $\bm{f}$, and then computing the Frobenius norm of the difference between the Jacobian of the sample and the nominal one, normalized with respect to the latter. This results into a set of \hl{indices} $\{\nu_{JS}^{(i)}\}$. The selected $\nu_{JS}$ is then the maximum among the $\{\nu_{JS}^{(i)}\}$. The goal of this section is to understand if our approach gives compatible results, and if it provides any advantage over the reference method.

The test case presented here is taken from~\citet{Junkins2004} and considers the propagation of an orbital uncertainty set in Cartesian coordinates along an eccentric orbit with semi-major axis  of 12871.5~km, eccentricity of 0.3629, and orbital inclination equal to 42.5~deg. An uncertainty of 1~km and 1~m/s in modulus for position and velocity is considered. These values are scaled down by a factor 10 with respect to the mentioned reference test to properly compare $\nu$ with $\nu_{JS}$. Indeed, the \acrshort{da}-based \acrshort{nli} aims at providing a local rigorous estimation of the nonlinearity of the transformation. This estimation is no longer accurate when the uncertainty set is large as terms with order higher than two cannot be neglected. This condition, however, does not represent a limitation for our approach since the automatic domain splitting routine is meant to prevent this from happening. As a result, if we want to have a direct comparison between $\nu$ and $\nu_{JS}$ over a single set, we need to select it such that the first is rigorous over its entire domain. A second aspect that is worth mentioning is the shape of the uncertainty region we are considering. As described in Section~\ref{sec:Index}, the set $[\bm{x}]$ is initialized by considering each component independently, and then building a symmetric interval around each nominal value. This can always be done via a priori eigendecomposition (see Section~\ref{sec:LOADS}). As a result, our uncertainty set is an $n$-dimensional hypercube around the nominal value of $\bm{x}$. If we want to sample the worst case surface to compute $\nu_{JS}$, we need to build a proper grid of points uniformly spaced over the $(n-1)$-dimensional surfaces of this hypercube.
\begin{figure}[!t]
\centering
\includegraphics[trim=2cm 0cm 2cm 1cm, clip=true, width=0.8\textwidth]{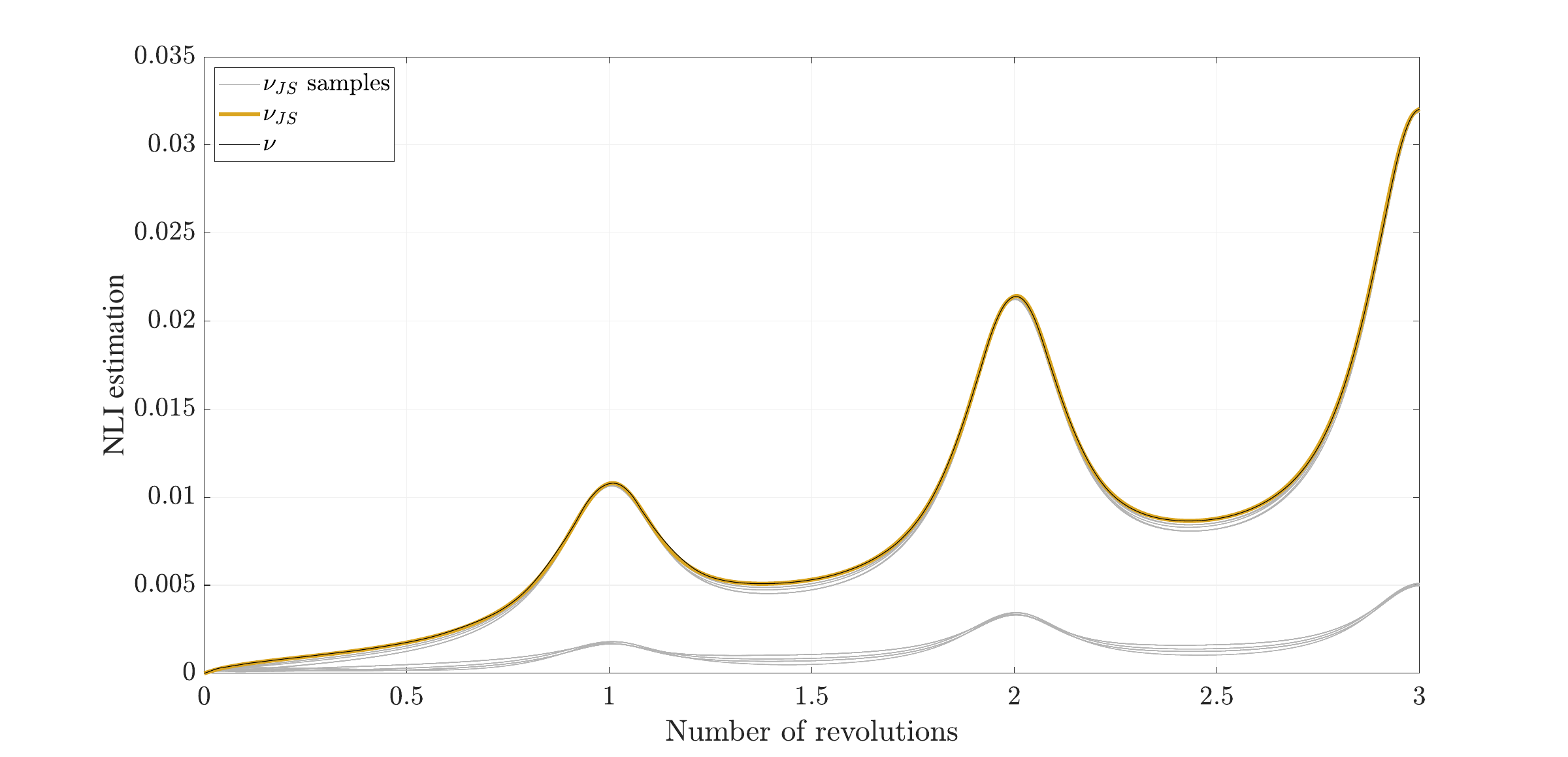}\\
\caption{Comparison between $\nu$ and $\nu_{JS}$ over three revolutions. \hl{The plot reports 50} $\nu_{JS}$ samples out of 150 for visualization purposes. }
\label{fig:NLI}
\end{figure}

Figure~\ref{fig:NLI} shows the result of the comparison. The plot shows the evolution of the \acrshortpl{nli} as a function of the number of revolutions. The $\nu_{JS}$ index converges with 150 samples uniformly distributed over the surface of the 6D uncertainty hypercube. The estimated $\nu_{JS}^{(i)}$ of the corresponding samples are plotted in grey, whereas the resulting \acrshort{nli} estimate is reported in gold. The trend of $\nu$ is instead reported in black. As can be seen, the matching between the two \hl{indices} is quite good throughout the propagation. This evidence is confirmed by Table~\ref{tab:NLI_trend}. The table illustrates the trend of the relative error between $\nu$ and $\nu_{JS}$ as a function of time. The error was obtained by considering 100 epochs per orbital period, computing the corresponding $\nu_j$ and $\nu_{JS,j}$ for each epoch $t_j$, and then computing the root mean square of the relative errors, i.e.
\begin{equation}
\textrm{RMSE}_{\nu}(t_j) = \sqrt{\dfrac{1}{n_j}\sum\limits_{j=1}^{n_j}\left(\dfrac{\nu_j-\nu_{JS,j}}{\nu_{JS,j}}\right)^2}
\end{equation}
where $n_j$ is the number of time epochs from 0 to $t_j$. As can be seen, the error is always lower than 4$\%$ and has larger fluctuations at the beginning, which progressively reduce as the propagation time increases. Overall, the matching between the two estimates is apparent. 

\hl{One} advantage of our approach is its computational efficiency. Regardless of the propagation epoch, the ratio between the runtime required by the \acrshort{da}-based propagation and a \hl{single} pointwise propagation of the orbital state oscillates between 2 and 3. This value can be interpreted as an equivalent number of samples $\nu_{JS}^{(i)}$ that could be \hl{computed} in the same amount of time required for the \hl{estimation} of $\nu$. Since convergence of $\nu_{JS}$ \hl{is} obtained with about 150 samples, we can say that our approach allows us to reduce the computational burden by at least a factor of \hl{50}. \hl{This factor is obviously problem dependent and could be significantly reduced by running the method by~\citet{Junkins2004} over multiple threads, which is not possible with our approach. However, it} comes hand in hand with the possibility of avoiding any sampling of the worst case surface, which makes our approach appealing for any generic application.

\begin{table}[!t]
\caption{Relative error between $\nu$ and $\nu_{JS}$ as a function of the propagation epoch.}
\label{tab:NLI_trend}
\begin{center}
\begin{tabular}{*{11}{c}}
\hline\noalign{\smallskip}
Epoch& 0.5 $T$& 1 $T$& 1.5 $T$& 2.00 $T$& 2.5 $T$& 3 $T$& 3.5 $T$& 4.00 $T$\\
\noalign{\smallskip}\hline\noalign{\smallskip}
$\textrm{RMSE}_{\nu}$ ($\%$)& 3.98& 2.91& 2.42& 2.10& 1.88& 1.72& 1.60& 1.52\\
\noalign{\smallskip}\hline
\end{tabular}
\end{center}
\end{table}

\subsection{The \acrshort{loads} algorithm for orbital uncertainty propagation}
\label{subsec:LOADS_test}
Once assessed the quality of the \acrshort{nli} estimation process, we illustrate here the application of the \acrshort{loads} algorithm to the case of orbital uncertainty propagation. We consider an equatorial orbit subject to drag and $J_2$ perturbations with a pericenter altitude of 300~km and an orbital eccentricity of 0.5. The object is at the pericenter at $t_0$. The dynamics are~\citep{DeMars2013,Gondelach2018}
\begin{equation}
\label{eq:EoM}
\bm{x}(t) = \begin{Bmatrix}
x\\y\\v_x\\v_y\\
\end{Bmatrix}
\qquad \dot{\bm{x}}(t) = \begin{Bmatrix}
v_x\\
v_y\\
-\mu_{\oplus}\dfrac{x}{r^3}-\dfrac{3J_2\mu_{\oplus} R^2_{\oplus}}{2 r^5}x-\dfrac{1}{2}\rho(h)\beta v_{rel}v_{rel,x}\\
-\mu_{\oplus}\dfrac{y}{r^3}-\dfrac{3J_2\mu_{\oplus} R^2_{\oplus}}{2 r^5}y-\dfrac{1}{2}\rho(h)\beta v_{rel}v_{rel,y}
\end{Bmatrix}
\end{equation}
where $x$ and $y$ are the two Cartesian coordinates describing the position vector in the Earth Centred Inertial (ECI) reference frame, $v_x$ and $v_y$ the velocity components, $r=\sqrt{x^2+y^2}$, $v_{rel,x}=v_x+\omega y$, $v_{rel,y}= v_y-\omega x$, $v_{rel}=\sqrt{v_{rel,x}^2+v_{rel,y}^2}$, $\omega$ is the angular velocity of the Earth, $\beta$ is the ballistic coefficient assumed equal to 1.4, whereas $\mu_{\oplus}$ and $R_{\oplus}$ are the Earth gravitational parameter and radius, respectively. The atmospheric density is described by an exponential model, where $\rho(h)=\rho_0\,\textrm{exp}\{(h_0-h)/h_s\}$, $\rho_0=3.614\cdot 10^{-14}$~kg/m$^3$, $h_0 = 700$~km and $h_s=88.667$~km. The data are taken from~\citet{DeMars2013}. The dynamics equations are written in Cartesian coordinates to magnify the nonlinearities, but the method can be  applied to any state representation.

\hl{The initial state is assumed to be normally distributed, with} uncertainties along the $x$ and $y$ directions only to help visualise the results. \hl{More specifically, }standard deviations of 10~km in $x$ and 100~km in $y$ are assumed. \hl{Such large values are selected to help visualize the performance of the \acrshort{loads} algorithm, and represent a realistic outcome of initial orbit determination with short and very short-arcs observations~\citep{Pirovano2021, Losacco2023}.} The \acrshort{loads} solver is then initialized by assuming a maximum number of splits per direction $N_{max,j}$ equal to 10 and a nonlinearity threshold $\varepsilon_{\nu}$ equal to 0.02. \hl{The merging algorithm is activated} with a number of merging breaks per revolution $n_{b/r}$ equal to 20. The uncertainty set is initialized with $\alpha_{\lambda_j} = 3$ and then propagated for two revolutions and a half.

\begin{figure}[!t]
\centering
\subfloat[\label{subfig:LOADS_t_33}]{\includegraphics[trim=2cm 0cm 2cm 1cm, clip=true, width=0.8\textwidth]{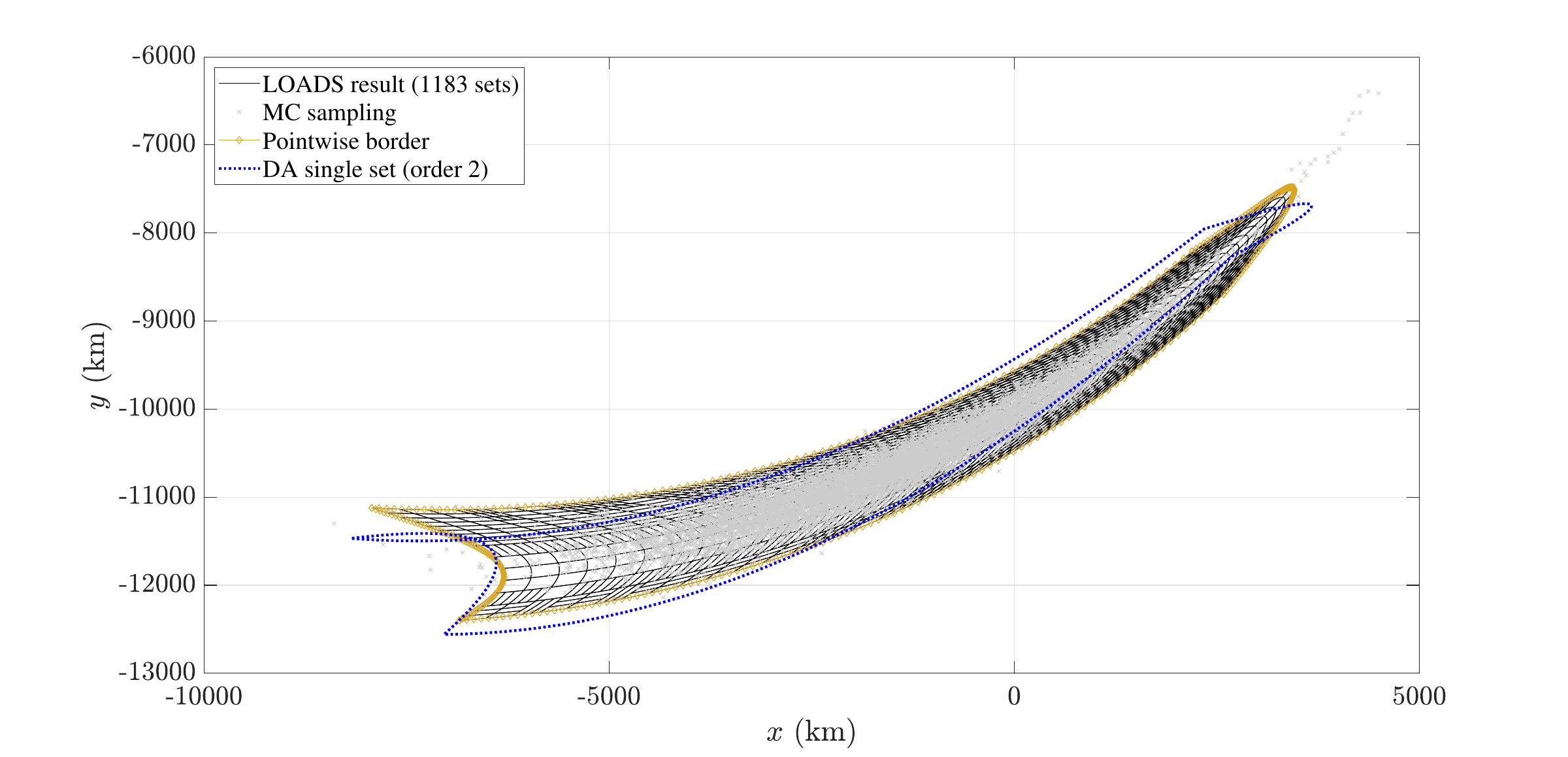}}\\
\subfloat[\label{subfig:LOADS_t_42}]{\includegraphics[trim=2cm 0cm 2cm 1cm, clip=true, width=0.8\textwidth]{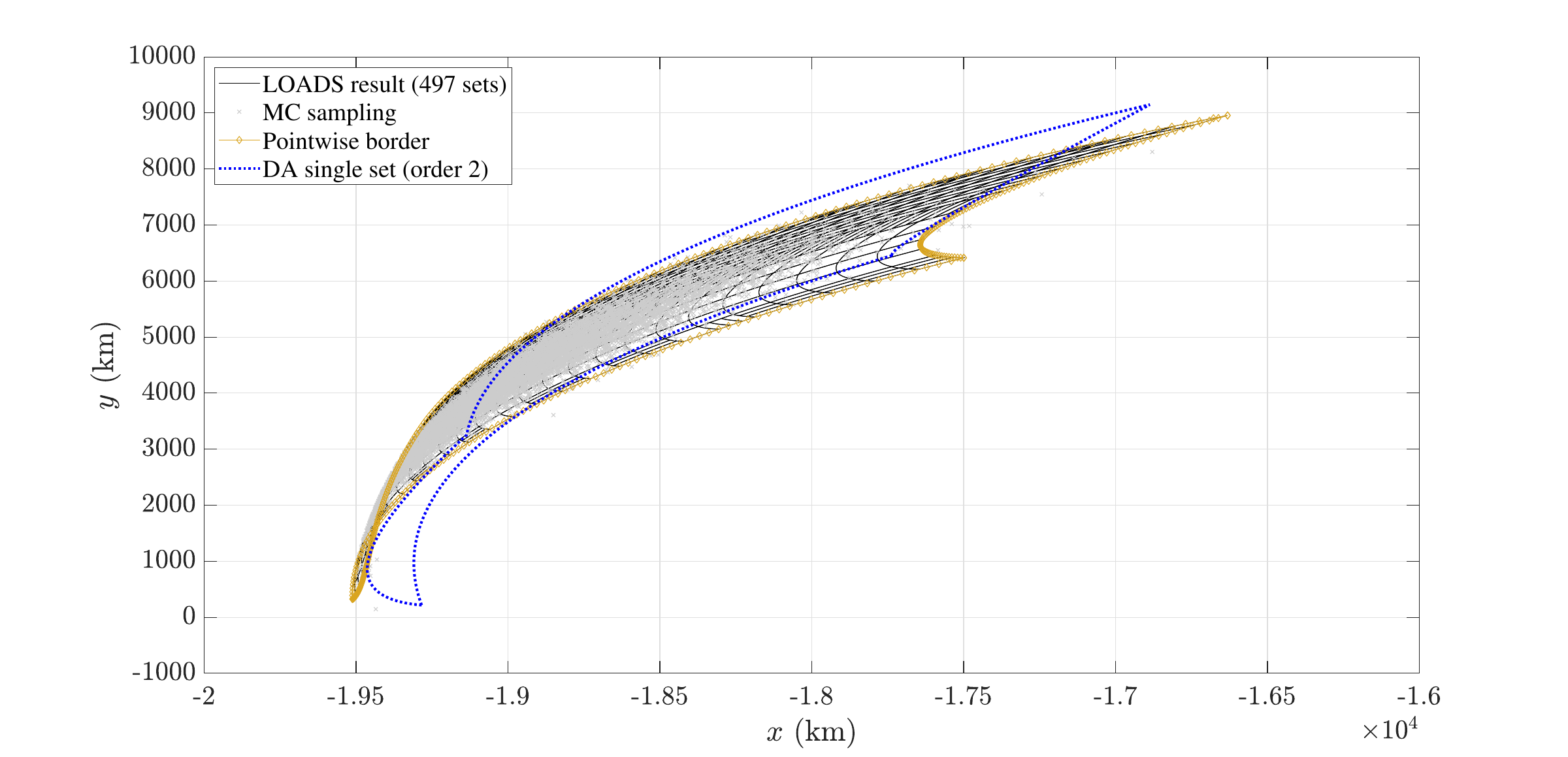}}\\
\caption{Orbital uncertainty propagation results: \protect\subref{subfig:LOADS_t_33} 1.88 orbital periods, pericenter approach, 1183 sets; \protect\subref{subfig:LOADS_t_42} 2.4 orbital periods, apocenter approach, 497 sets ($\varepsilon_{\nu}=0.02$, $N_{max,j}=10$, $\alpha_{\lambda_j}=3$, $n_{b/r}=20$).}
\label{fig:LOADS}
\end{figure}
\begin{figure}[!t]
\centering
\subfloat[\label{subfig:LOADS_t_33_focus}]{\includegraphics[trim=2cm 0cm 2cm 0.5cm, clip=true, width=0.8\textwidth]{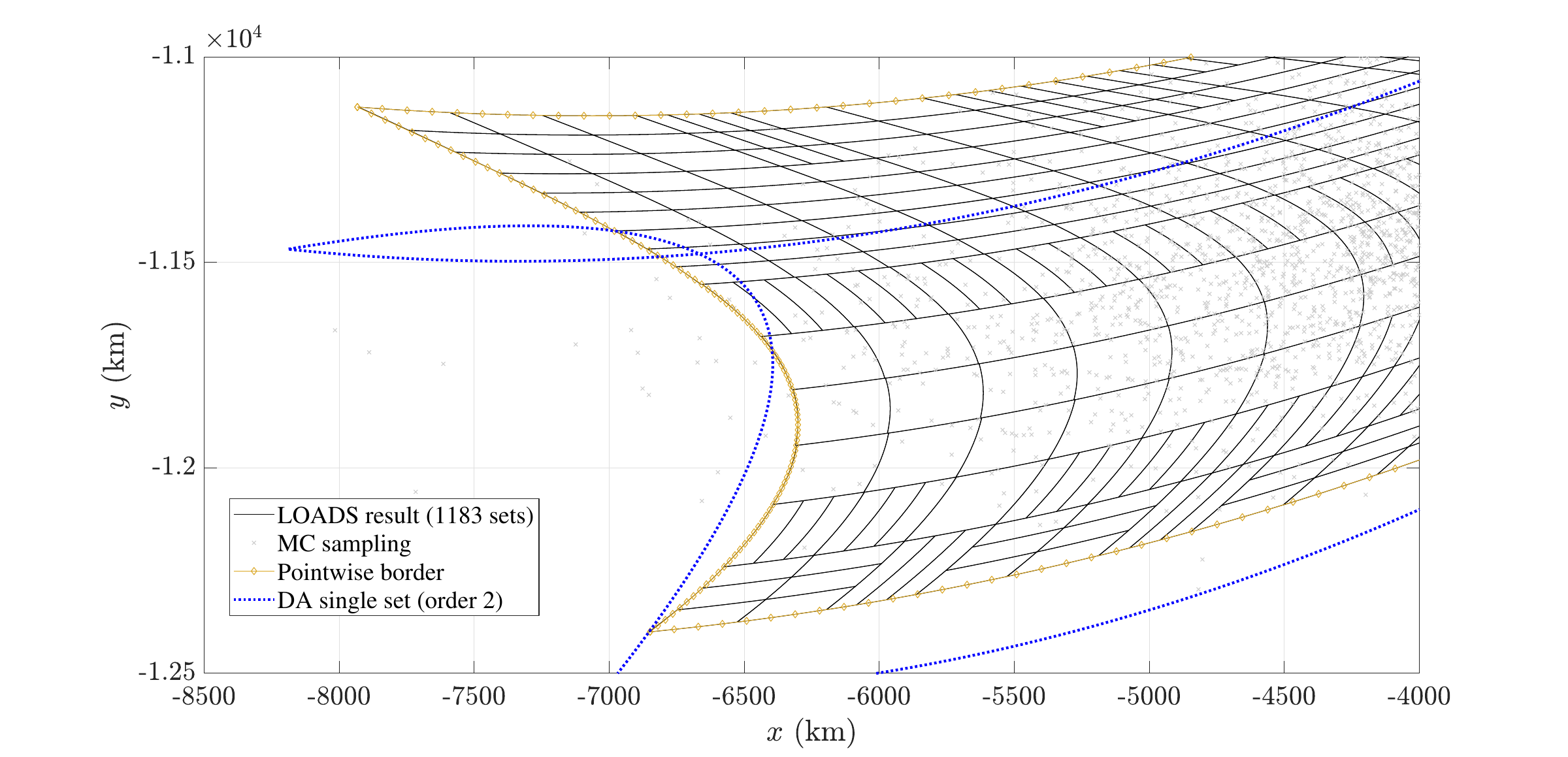}}\\
\subfloat[\label{subfig:LOADS_t_42_focus}]{\includegraphics[trim=2cm 0cm 2cm 1cm, clip=true, width=0.8\textwidth]{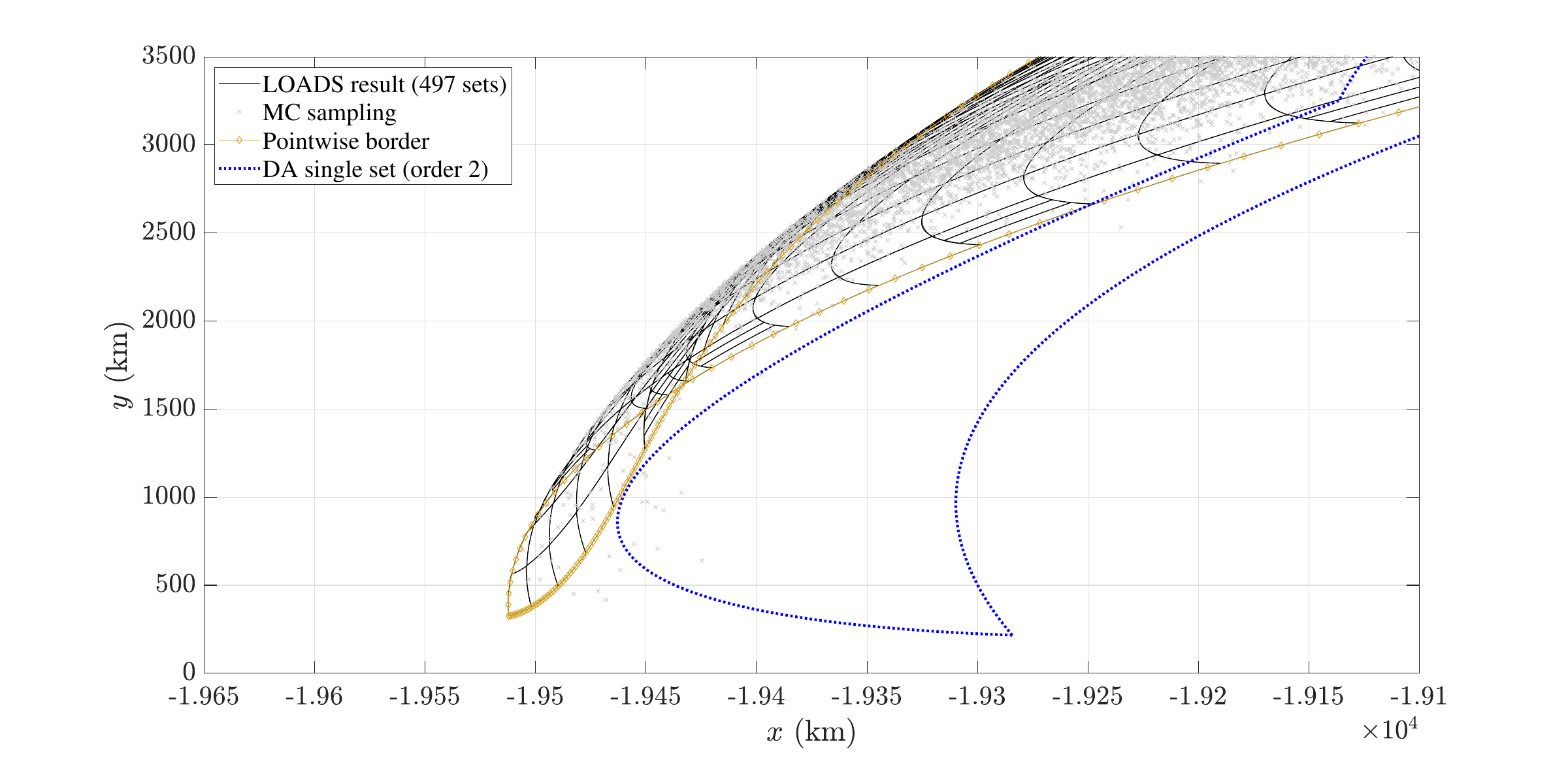}}\\
\caption{Orbital uncertainty propagation results (detail): \protect\subref{subfig:LOADS_t_33_focus} 1.88 orbital periods, pericenter approach, 1183 sets; \protect\subref{subfig:LOADS_t_42_focus} 2.4 orbital periods, apocenter approach, 497 sets ($\varepsilon_{\nu}=0.02$, $N_{max,j}=10$, $\alpha_{\lambda_j}=3$, $n_{b/r}=20$).}
\label{fig:LOADS_focus}
\end{figure}

The results of the simulation are illustrated in Fig.~\ref{fig:LOADS}. The plot shows the evolution of the uncertainty set at two epochs: after 1.88 revolutions (1.88 $T$ , Fig.~\ref{subfig:LOADS_t_33}, pericenter approach) and 2.4 revolutions (2.4 $T$, Fig.~\ref{subfig:LOADS_t_42}, apocenter approach). For each plot, the results of four different simulations are reproduced. The gold diamonds describe a pointwise propagation of the borders of the $3\sigma$ initial uncertainty set\hl{, i.e. the result of the one-by-one propagation from $t_0$ to $t_f$ of samples uniformly spaced over the $3\sigma$ enclosure of the initial uncertainty set}. The grey dots represent the projection onto the $x-y$ plane of \acrshort{mc} samples drawn from the initial multivariate normal distribution and propagated to the final epoch. The blue dotted border is the result of a single-set second-order propagation of the $3\sigma$ uncertainty box. In contrast, the black grid represents the sets generated by the \acrshort{loads} algorithm when propagating and projecting the same initial box. That is, the grey dots and gold diamonds represent the real orbital uncertainty distribution and 3$\sigma$ enclosure, respectively, whereas the black grid and the blue box the attempt to approximate \hl{this enclosure with second-order Taylor polynomials}, with and without domain splitting. 

Let us first compare the actual 3$\sigma$ enclosure (gold diamonds) with its second-order, single-set \hl{approximation} (dotted blue line). As can be seen, a single \hl{polynomial} fails to capture the exact curvature of the uncertainty set. This feature is particularly evident at the \hl{extrema}, where the propagated uncertainty set significantly bends, and its single-set \hl{approximation} diverges. This behaviour becomes more apparent for longer propagation windows, as can be seen by comparing Fig.~\ref{subfig:LOADS_t_33} and Fig.~\ref{subfig:LOADS_t_42}. The \acrshort{loads} algorithm aims at maintaining the computational agility of second-order expansions while increasing the overall accuracy. This is achieved by approximating the uncertainty set with a manifold of locally second-order expansions. The comparison between the blue dotted line and the black grid clearly shows the capability of the \acrshort{loads} to accurately describe the curvature of the set using multiple domains. The local shape of the single boxes can be better inspected by looking at Fig.~\ref{subfig:LOADS_t_33_focus}, which shows a detail of the most stretched portion of Fig.~\ref{subfig:LOADS_t_33}. The generated sets are quasi-parallelograms, with smaller size where the curvature of the overall set increases. \hl{The propagation process causes the uncertainty set to stretch and eventually fold over itself, as can be seen by inspecting the bottom-left portion of Fig.~\ref{subfig:LOADS_t_42_focus}. Still, the \acrshort{loads} can accurately describe the shape of the propagated uncertainty set. Note also that some \acrshort{mc} samples actually fall out of the propagated set. This can be seen both in Fig.~\ref{subfig:LOADS_t_33_focus} ($x$<-6500~km, -1.2e4~km<$y$<-1.1e4~km) and in Fig.~\ref{subfig:LOADS_t_42_focus} (-1.95e4~km<$x$<-1.94e4~km, 500~km<$y$<1000~km). This is not unexpected, as the selected $3\sigma$ enclosure unavoidably discards a minimal but non irrelevant fraction of the probability mass.}

The accuracy of the representation is qualitatively confirmed by the distribution of the \acrshort{mc} samples and quantified in Table~\ref{tab:LOADS_vs_MC}. The table shows a comparison between the single-set (\textit{``1 set''} superscript) approach and the \acrshort{loads} algorithm in terms of errors in position $\varepsilon_{p}^{bd}$ and velocity $\varepsilon_{v}^{bd}$. These errors are computed by uniformly sampling the edges of each set (one for the first approach, several for the \acrshort{loads}) with $n_{edge}$ samples per edge, and then computing the root mean square of the difference between pointwise propagation (our reference) and polynomial evaluation. More specifically, let us define with $\bm{x}^{(k)}_{0,l}$ the generic sample drawn at time $t_0$ from the set $k$ of $M_0$ and $\left[\bm{x}^{(k)}(t_f)\right]$ the element of $M_f$ mapping the sample from $t_0$ to $t_f$. Let us define with $\hat{\bm{x}}_{f,l}^{(k)}$ the mapping of $\bm{x}^{(k)}_{0,l}$ through $\left[\bm{x}^{(k)}(t_f)\right]$ and $\bm{x}_{f,l}^{(k)}$ the result of pointwise propagation, i.e. $\bm{x}_{f,l}^{(k)}=\bm{G}\left(\bm{x}^{(k)}_{0,l},t_f\right)$. Let us now consider the generic root mean square error as
\begin{equation}
\varepsilon_{\bm{x}}^{bd}=\left\{\varepsilon_{x}^{bd},\varepsilon_{y}^{bd},\varepsilon_{v_x}^{bd},\varepsilon_{v_y}^{bd}\right\}^{\textrm{T}}=\sqrt{\dfrac{1}{n_{set}n_{edge}}\sum\limits_{k=1}^{n_{set}}\sum\limits_{l=1}^{n_{edge}}\left(\hat{\bm{x}}_{f,l}^{(k)}-\bm{x}_{f,l}^{(k)}\right)\odot\left(\hat{\bm{x}}_{f,l}^{(k)}-\bm{x}_{f,l}^{(k)}\right)}
\end{equation}
where $\odot$ indicates the Hadamar\hl{d} (component-wise) product between the two vectors, while the double summation is performed element-wise, and $n_{set}$ is the number of sets. As a result, the errors $\varepsilon_{p}^{bd}$ and $\varepsilon_{v}^{bd}$ are computed as
\begin{equation}
\begin{gathered}
\varepsilon_{p}^{bd} = \textrm{max}\left\{\varepsilon_{x}^{bd},\varepsilon_{y}^{bd}\right\}\\
\varepsilon_{v}^{bd} = \textrm{max}\left\{\varepsilon_{v_x}^{bd},\varepsilon_{v_y}^{bd}\right\}\\
\end{gathered}
\end{equation}

\begin{table}[!t]
\caption{\Acrshort{loads} and single set accuracies when compared against pointwise propagation of samples drawn at the borders of the generated sets ($\varepsilon_{\nu}=0.02$, $N_{max,j}=10$, $\alpha_{\lambda_j}=3$, $n_{b/r}=20$).}
\label{tab:LOADS_vs_MC}
\begin{center}
\begin{tabular}{*{6}{c}}
\hline\noalign{\smallskip}
\begin{tabular}{@{}c@{}}Epoch \\ (T)\end{tabular}&
\begin{tabular}{@{}c@{}}$\varepsilon_{p}^{bd,\textit{1}\,set}$\\ (km)\\\end{tabular}& 
\begin{tabular}{@{}c@{}}$\varepsilon_{v}^{bd,\textit{1}\,set}$\\ (m/s)\\\end{tabular}& 
\begin{tabular}{@{}c@{}}$n_{set}^{LOADS}$\\ \\\end{tabular}&
\begin{tabular}{@{}c@{}}$\varepsilon_{p}^{bd,LOADS}$\\ (km)\\\end{tabular}& 
\begin{tabular}{@{}c@{}}$\varepsilon_{v}^{bd,LOADS}$\\ (m/s)\\\end{tabular}\\
\noalign{\smallskip}\hline\noalign{\smallskip}
0.88&     \phantom{0}52.541&           \phantom{0}35.695&  \phantom{0}315&  0.976&      0.493\\
1.40&      \phantom{0}39.450&  \phantom{0}\phantom{0}7.876&   \phantom{0}335& 0.668&     0.214\\
1.88&                 190.338&   166.368&    1183&  1.862&                              1.124\\
2.40&              117.327&   \phantom{0}23.094& \phantom{0}497&  1.143&                0.378\\
\noalign{\smallskip}\hline
\end{tabular}
\end{center}
\end{table}
The results of the single-set representation show two evident trends. The first one is that the accuracy decreases with the number of revolutions.  This drop in accuracy is apparent when comparing two consecutive passages approaching the pericenter (lines 1 and 3) or the apocenter (lines 2 and 4).  This result is expected as the uncertainty set tends to stretch, thus preventing the single second-order set from describing it accurately. The trend of the error, however, is not monotonic. Indeed, errors drastically decrease when passing from the pericenter to the apocenter (compare line 1 with line 2, or line 3 with line 4). This is an expected result as the nonlinearities are stronger close to the pericenter, and decrease while departing from it. 

All these elements directly affect the performance of the \acrshort{loads} algorithm. Three aspects can be noticed. The first one is a relevant improvement in the accuracy. By splitting the uncertainty set into smaller domains, the algorithm can describe it more accurately. Moreover, unlike the single set case, the obtained accuracy is almost unaffected by the propagation time. This property is enabled by the splitting algorithm, which automatically regulates the number of sets to meet the $\varepsilon_{\nu}$ constraint. The final trend can be noticed by looking at the number of sets. 
While is it almost unaltered when passing from 0.88 to 1.40, it significantly drops one revolution later, from 1.88 to 2.40. The trend is compliant with the ``accordion-like'' behaviour of the error previously mentioned, and it is well captured by the merging algorithm, which recombines the sets when possible.

A more detailed analysis of the advantages introduced by the merging breaks is offered in Table.~\ref{tab:Sets_recovery}. The table shows the evolution of the number of sets \hl{$n_{set}$ and the required computational time $t_{CPU}$} as a function of the propagation epoch (expressed in fractions of the period $T$) for two different approaches: a simple \acrshort{loads}, where no merging is performed, and the combined \acrshort{loads}-merging approach. \hl{Let us first analyse the number of sets.} When no merging is allowed, the number of sets continuously grows, with a steep growth at pericenter passages, and more moderate rates while departing from it. When merging breaks are introduced, recombination of many of the pericenter passage-induced sets occurs, reducing the total number of sets. New splits are performed when passing closer to the pericenter. However, the number of sets obtained with the \acrshort{loads}-merging algorithm at the pericenter passage is lower than what is obtained with the simple \acrshort{loads} approach. \hl{Overall, }the introduction of the merging breaks both guarantees a more efficient description of the uncertainty when nonlinearities are smaller and reduces the overall number of sets while approaching regions of larger nonlinearities.  

\begin{table*}[!t]
\caption{Number of generated sets \hl{$n_{set}$ and computational time $t_{CPU}$} as a function of the propagation epoch: no merging vs merging approach ($\varepsilon_{\nu}=0.02$, $N_{max,j}=10$, $\alpha_{\lambda_j}=3$, $n_{b/r}=100$). \hl{In case of merging, the table reports the impact of the merging phase in terms of percentage of the overall computational load $t_{CPU,m}$.}}
    \label{tab:Sets_recovery}
    \centering
    \small
    \sisetup{round-mode=places,round-precision=4}
    \begin{tabular}{c c c c c c c c c c c}
    \hline\hline
    & & 0.25 $T$& 0.5 $T$& 0.75 $T$& 1 $T$& 1.25 $T$& 1.5 $T$& 1.75 $T$& 2.00 $T$& 2.50 $T$\\
    \cline{3-11}
\multirow{2}{*}{\shortstack[c]{$n_{set}$}}&	No merging& 27& 59& 121& 1379& 1520& 1911& 2101& 5023& 8937\\
& Merging&  27& 45& \phantom{0}99& 1375& \phantom{0}343& \phantom{0}339& \phantom{0}395& 3597& \phantom{0}485\\
\hline
\multirow{3}{*}{\shortstack[c]{\hl{$t_{CPU}$ (s)}\\\hl{($t_{CPU,m}$ ($\%$))}}}&	\hl{No merging}& \hl{0.60}& \hl{1.44}& \hl{2.76}& \hl{9.58}& \hl{25.88}& \hl{43.54}& \hl{65.15}& \hl{99.73}& \hl{250.28}\\
& \multirow{2}{*}{\shortstack[c]{\hl{Merging}}}&  \multirow{2}{*}{\shortstack[c]{\hl{0.65}\\\hl{(1.56$\%$)}}}& \multirow{2}{*}{\shortstack[c]{\hl{1.33}\\ \hl{(1.62$\%$)}}}& \multirow{2}{*}{\shortstack[c]{\hl{2.71}\\\hl{(1.28$\%$)}}}& \multirow{2}{*}{\shortstack[c]{\hl{8.87}\\\hl{(0.96 $\%$)}}}& \multirow{2}{*}{\shortstack[c]{\hl{16.33}\\\hl{(0.89$\%$)}}}& \multirow{2}{*}{\shortstack[c]{\hl{20.06}\\\hl{(0.85$\%$)}}}& \multirow{2}{*}{\shortstack[c]{\hl{24.15}\\\hl{(0.83 $\%$)}}}& \multirow{2}{*}{\shortstack[c]{\hl{42.87}\\\hl{(0.76$\%$)}}}& \multirow{2}{*}{\shortstack[c]{\phantom{0}\hl{69.92}\\\hl{(0.74$\%$)}}}\\\\
\hline\hline
    \end{tabular}
    \vspace{0.2cm}
\end{table*}

\hl{The merging process not only simplifies the structure of the resulting manifold, but also speeds up the propagation process. Evidence for this is found by comparing the computational time required by the two approaches. As can be seen, when merging is active, the propagation is almost always faster, and the advantages grow as the time frame increases. This trend can be explained by looking at the relative impact of the merging phases on the overall computational load ($t_{CPU,m}$ parameter). The time required by the merging process ranges from 0.7$\%$ to 1.6$\%$ of $t_{CPU}$, with lower values corresponding to longer propagation windows, which require a larger number of sets and for which, therefore, the splitting process becomes dominant. The introduction of merging phases partially mitigates the growth of the sets in time with a limited effect on the overall computational time, thus resulting into the highlighted trend.}

\subsubsection{Selecting the \acrshort{loads} control parameters}
\label{subsub:LOADS_param}
This section investigates the impact of the nonlinearity threshold $\varepsilon_{\nu}$ and the size of the uncertainty domain on the \acrshort{loads} performance. Their role is studied in terms of both accuracy and  computational load. The quality of the results is studied by investigating the obtained statistical representation of the propagated uncertainty set. The analysis is done by drawing a set of $N_{MC}$ \acrshort{mc} samples $\left\{\bm{x}_{0,l}\right\}$ from the initial Gaussian distribution. Following a procedure similar to the one described in Section~\ref{subsec:LOADS_test}, each sample is then mapped and propagated to $t_f$, thus obtaining $\hat{\bm{x}}_{f,l}$ and $\bm{x}_{f,l}$, respectively. The comparison between the two sets of samples allows us to characterise the quality of \hl{the statistical} moments of the propagated \acrshort{pdf}. Six different additional \hl{indices} are considered. The first two \hl{indices} are $\varepsilon_{p}$ and $\varepsilon_{v}$, i.e. the root mean square errors in position and velocity between \acrshort{mc} and \acrshort{loads} samples. These two \hl{indices} do not coincide with $\varepsilon_{p}^{bd}$ and $\varepsilon_{v}^{bd}$, which were built by considering regular grids only along the borders. The third index, here referred as $I_{\mu}$, is the relative error in the norm of the \acrshort{pdf} mean $\bm{\mu}$, i.e.
\begin{equation}
I_{\mu} = \dfrac{\left\lvert\left\lVert\bm{\mu}^{LOADS}\right\rVert-\left\lVert\bm{\mu}^{MC}\right\rVert\right\rvert}{\left\lVert\bm{\mu}^{MC}\right\rVert}
\end{equation}
where $\bm{\mu}^{LOADS}$ and $\bm{\mu}^{MC}$ are the \acrshort{pdf} means computed with \acrshort{loads} and \acrshort{mc} samples, respectively. The fourth index is the relative error in the maximum eigenvalue of the sample covariance matrices
\begin{equation}
I_{\lambda} = \dfrac{\left\lvert\lambda_{max}^{LOADS}-\lambda^{MC}_{max}\right\rvert}{\lambda_{max}^{MC}}
\end{equation}
The fifth and six \hl{indices} are the relative errors in Mardia's multivariate skewness $\gamma_{1,p}$ and kurtosis $\gamma_{2,p}$, respectively
\hl{
\begin{subequations}
\begin{gather}
    I_{\gamma_{1,p}} = \dfrac{\left\lvert\gamma_{1,p}^{LOADS}-\gamma_{1,p}^{MC}\right\rvert}{\gamma_{1,p}^{MC}}\\
    I_{\gamma_{2,p}} = \dfrac{\left\lvert\gamma_{2,p}^{LOADS}-\gamma_{2,p}^{MC}\right\rvert}{\gamma_{2,p}^{MC}}
\end{gather}    
\end{subequations}
}
where the generic sample skewness and kurtosis are computed as
\hl{
\begin{subequations}
\begin{gather}
    \gamma_{1,p} = \dfrac{1}{N_{MC}^2}\displaystyle\sum\limits_{q=1}^{N_{MC}}\sum\limits_{r=1}^{N_{MC}}\left[\left(\bm{x}_q-\bm{\mu}\right)^{\textrm{T}}{\bm{P}^{-1}}\left(\bm{x}_r-\bm{\mu}\right)\right]^3\\
    \gamma_{2,p} = \dfrac{1}{N_{MC}}\displaystyle\sum\limits_{q=1}^{N_{MC}}\left[\left(\bm{x}_q-\bm{\mu}\right)^{\textrm{T}}{\bm{P}^{-1}}\left(\bm{x}_q-\bm{\mu}\right)\right]^2
\end{gather}
\end{subequations}
}
with $\bm{P}$ representing the sample covariance matrices~\citep{Mardia1970, Yuan2004}.
In all the presented analyses, the number of \acrshort{mc} samples (\num{2e4}) was selected as the one granting the statistical convergence \hl{(second decimal place)} of the considered \hl{indices}.

\begin{table}[!t]
\caption{\acrshort{loads} performance as a function of $\varepsilon_{\nu}$ (1.88 revolutions, $N_{max,j}=10$, $\alpha_{\lambda_j}=3$, $n_{b/r}=20$).}
\label{tab:LOADS_vs_MC_nu}
\begin{center}
\begin{tabular}{*{11}{c}}
\hline\noalign{\smallskip}
\begin{tabular}{@{}c@{}}$\varepsilon_{\nu}$\\ \\\end{tabular}&
\begin{tabular}{@{}c@{}}$n_{set}^{LOADS}$\\ \\\end{tabular}&
\begin{tabular}{@{}c@{}}$\tilde{t}_{CPU}^{LOADS}$\\ \\\end{tabular}&
\begin{tabular}{@{}c@{}}$\varepsilon_{p}^{bd}$\\ (km)\\\end{tabular}&
\begin{tabular}{@{}c@{}}$\varepsilon_{v}^{bd}$\\ (m/s)\\\end{tabular}&
\begin{tabular}{@{}c@{}}$\varepsilon_{p}$\\ (km)\\\end{tabular}&
\begin{tabular}{@{}c@{}}$\varepsilon_{v}$\\ (m/s)\\\end{tabular}&
\begin{tabular}{@{}c@{}}$I_{\mu}$\\ ($\%$)\\\end{tabular}&
\begin{tabular}{@{}c@{}}$I_{\lambda}$\\ ($\%$)\\\end{tabular}&
\begin{tabular}{@{}c@{}}$I_{\gamma_{1,p}}$\\ ($\%$)\\\end{tabular}&
\begin{tabular}{@{}c@{}}$I_{\gamma_{2,p}}$\\ ($\%$)\\\end{tabular}\\
\noalign{\smallskip}\hline\noalign{\smallskip}
0.02& 1183& 1.000& \phantom{0}\phantom{0}1.86& \phantom{0}\phantom{0}1.12& \phantom{0}2.23& \phantom{0}1.24& 2.63e-4& 1.28e-2& \phantom{0}1.01& \phantom{0}1.49\\
0.04& \phantom{0}507& 0.377& \phantom{0}19.81& \phantom{0}10.85& \phantom{0}9.89& \phantom{0}5.43& 1.85e-3& 1.70e-2& \phantom{0}3.52& \phantom{0}5.12\\
0.06& \phantom{0}427& 0.330& \phantom{0}47.04& \phantom{0}28.11& 23.41& 13.84& 1.17e-2& 2.44e-1& \phantom{0}9.85& 14.33\\
$\infty$& \phantom{0}\phantom{0}\phantom{0}1& 0.008& 208.53& 160.26& 31.47& 39.991& 1.44e-2& 1.72e0& 56.72& 59.55\\
\noalign{\smallskip}\hline
\end{tabular}
\end{center}
\end{table}

The first parameter here investigated is the \acrshort{nli} threshold $\varepsilon_{\nu}$. As described in Section~\ref{sec:LOADS}, this parameter provides an upper bounds for $\nu$: when $\nu>\varepsilon_{\nu}$, splits are required. As a result, a finer mesh should be obtained by lowering the threshold, which should imply a higher accuracy in the obtained results. An analysis of the impact of the nonlinearity threshold on the accuracy of the results is shown in Table~\ref{tab:LOADS_vs_MC_nu}. The table shows how the \acrshort{loads} performance changes for three values of $\varepsilon_{\nu}$: 0.02 (our nominal solution), 0.04 and 0.06. In addition, the results obtained with a single second-order set ($\varepsilon_{\nu}=\infty$) are presented. The performance is expressed in terms of the six described \hl{indices}. In addition, the number of sets $n_{set}^{LOADS}$, the errors in position and velocity along the borders $\varepsilon_{p}^{bd}$ and $\varepsilon_{v}^{bd}$, and the required computational time $\tilde{t}_{CPU}^{LOADS}$, normalized with respect to the $\varepsilon_{\nu}= 0.02$ case, are shown. 

\begin{figure}[!t]
\centering
\subfloat[\label{subfig:LOADS_t_33_2sigma}]{\includegraphics[trim=2cm 0cm 2cm 1cm, clip=true, width=0.8\textwidth]{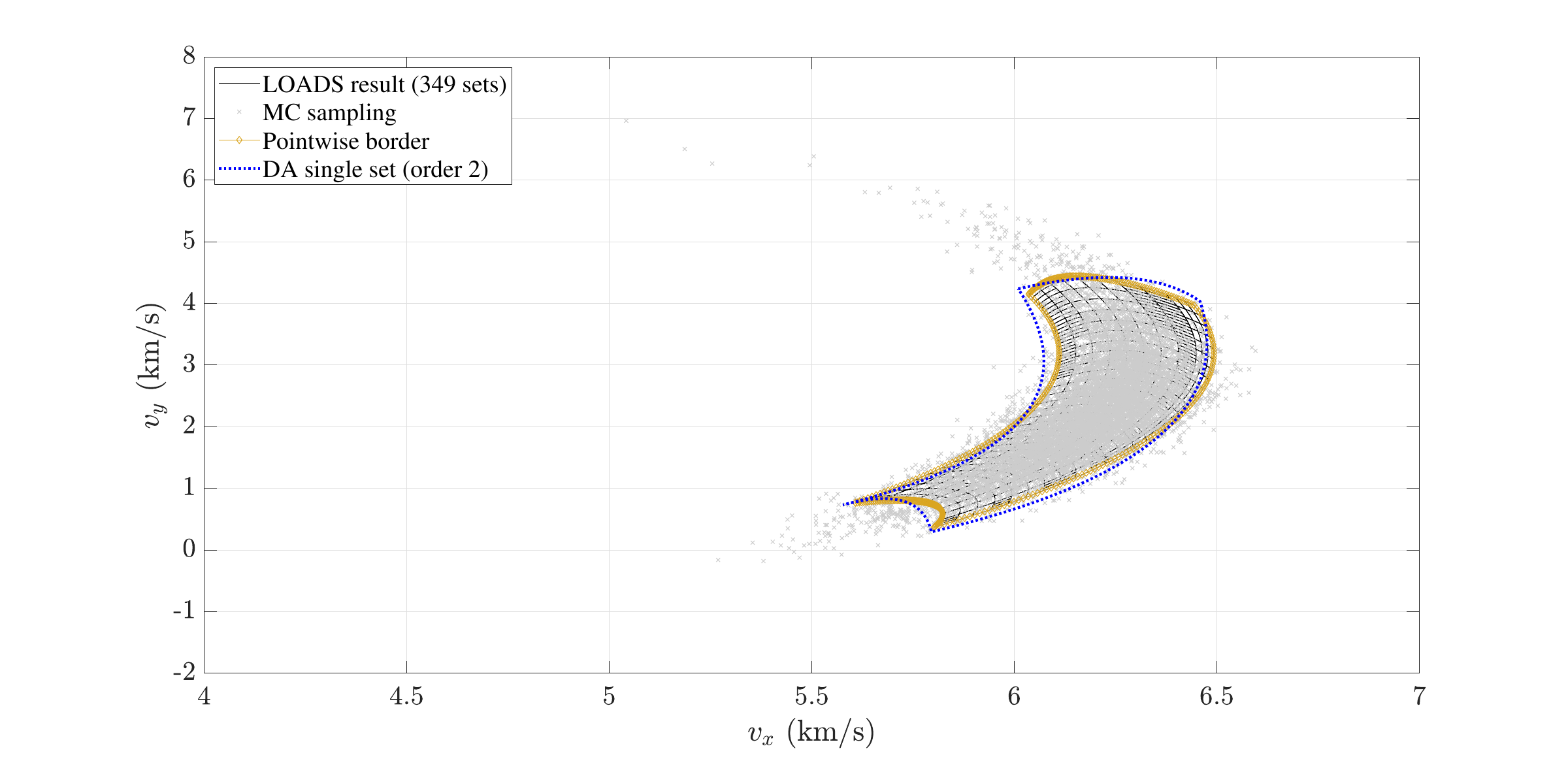}}\\
\subfloat[\label{subfig:LOADS_t_33_4sigma}]{\includegraphics[trim=2cm 0cm 2cm 1cm, clip=true, width=0.8\textwidth]{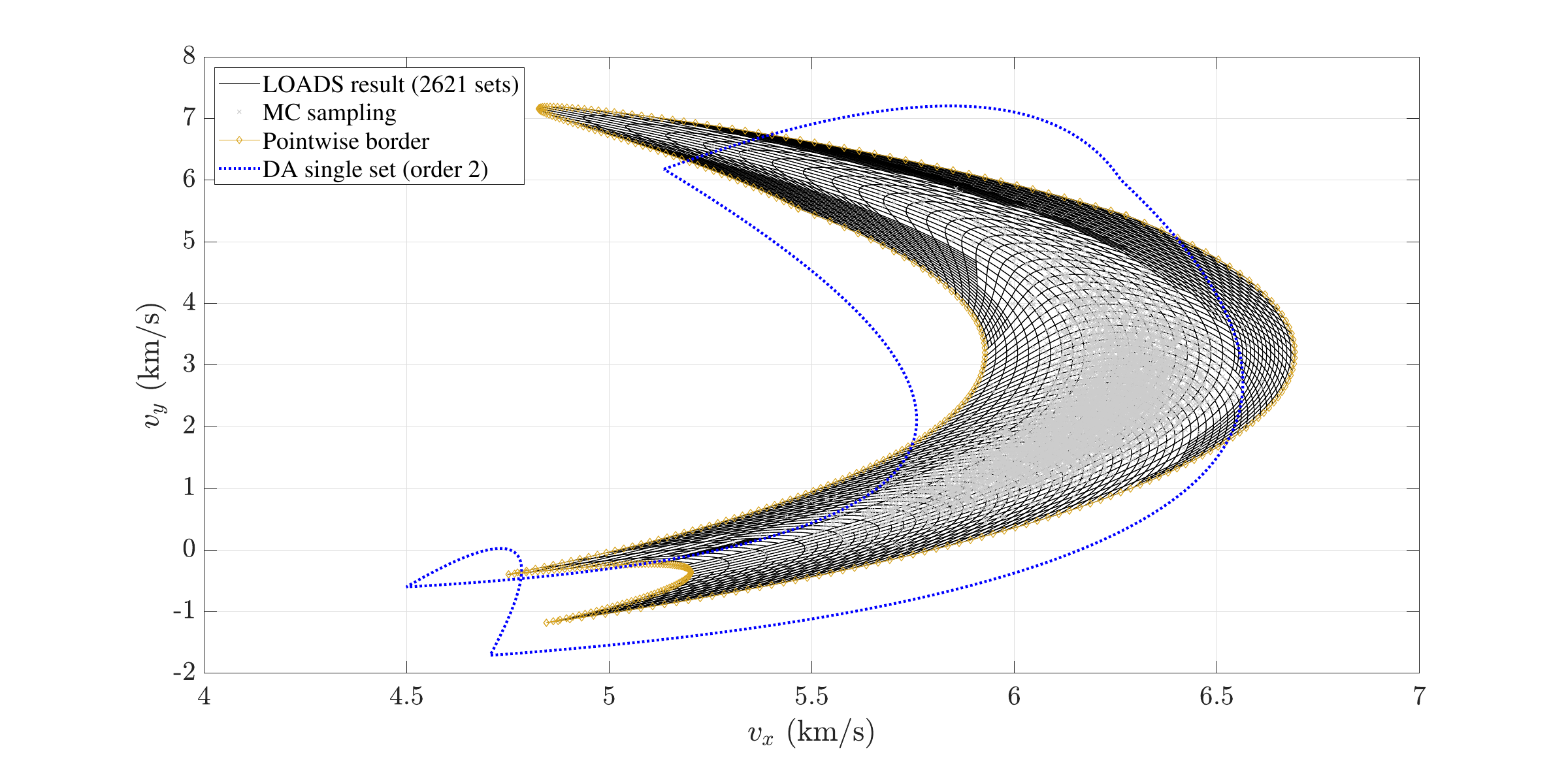}}\\
\caption{Orbital uncertainty propagation results after 1.88 orbital revolutions: \protect\subref{subfig:LOADS_t_33_2sigma} 2$\sigma$; \protect\subref{subfig:LOADS_t_42} $4\sigma$ domain ($\varepsilon_{\nu}=0.02$, $N_{max,j}=10$, $n_{b/r}=20$).}
\label{fig:LOADS_sigma}
\end{figure}
As can be seen, by increasing the value of the threshold $\varepsilon_{\nu}$, the number of generated sets and the required computational time drastically decrease, with a reduction of around a factor 3 in both \hl{indices} while passing from $\varepsilon_{\nu}=0.02$ to 0.06. This, however, comes hand in hand with a reduction in the obtained accuracy. More specifically, the errors in position and velocity along the borders increase of more than a factor 20, while an increase of about one order of magnitude in the $\varepsilon_{p}$ and $\varepsilon_{v}$ \hl{indices} occurs. These two \hl{indices} are however generally lower than $\varepsilon_{p}^{bd}$ and $\varepsilon_{v}^{bd}$, as the latter include more samples in the extreme regions of the set. Similar trends can be identified for the $I_{\mu}$, $I_{\lambda}$, $I_{\gamma_{1,p}}$ and $I_{\gamma_{2,p}}$ \hl{indices}, which significantly increase with the value of $\varepsilon_{\nu}$. The last line of the table illustrates the performance in case no split is performed, with huge errors in all the parameters. This is of course expected, as a single polynomial expansion cannot accurately cover the evolution of the uncertainty set. Overall, if an accurate description of both low and high order moments of the propagated \acrshort{pdf} is desired, the value of $\varepsilon_{\nu}$ shall be maintained low. The selection of $\varepsilon_{\nu}$ is application dependent and can be made less strict by our use of second-order expansions, which allows us to better tolerate mild nonlinearities.

The second parameter affecting the accuracy of the results is the size of the uncertainty set. The \acrshort{loads} algorithm, indeed, is initialized by defining a specific scaling factor $\alpha_{\lambda_j}\sqrt{\lambda_j}$ for each variable. Then, the bounds for the Jacobian matrix are computed by assuming that each variable can vary in the range $[-1,+1]$. By increasing the value of $\alpha_{\lambda_j}$, the amount of probability mass included in the initial uncertainty set increases. Thus, a better description of the propagated \acrshort{pdf} should be obtained. Conversely, a larger size should result \hl{in} a greater number of sets.

The first confirmation of this trend can be found in Fig.~\ref{fig:LOADS_sigma}. The two plots show the results of the \acrshort{loads} algorithm after 1.88 revolutions when projected onto the $v_x-v_y$ plane considering the nominal control parameters ($\varepsilon_{\nu}=0.02$, $N_{max,j}=10$, $n_{b/r}=20$), and a different size of the uncertainty set: $2\sigma$ (Fig.~\ref{subfig:LOADS_t_33_2sigma}) and $4\sigma$ (Fig.~\ref{subfig:LOADS_t_33_4sigma}). As can be seen, when a $2\sigma$ uncertainty set is considered, a smaller area is covered. Thus, a low number of sets (349) is generated, with the single set result reproducing the shape of the propagated uncertainty quite well. However, a large portion of the generated \acrshort{mc} samples falls out of the considered domain. These samples are instead included in the $4\sigma$ region, which requires a much larger number of second-order sets to be accurately described (2621).
\begin{table}[!t]
\caption{\Acrshort{loads} performance as a function of the domain size $\alpha_{\lambda_j}$ ($\varepsilon_{\nu}=0.02$, $N_{max,j}=20$, $n_{b/r}=20$).}
\label{tab:LOADS_vs_MC_sigma}
\begin{center}
\begin{tabular}{*{11}{c}}
\hline\noalign{\smallskip}
\begin{tabular}{@{}c@{}}$\alpha_{\lambda_j}$\\ \\\end{tabular}&
\begin{tabular}{@{}c@{}}$n_{set}^{LOADS}$\\ \\\end{tabular}&
\begin{tabular}{@{}c@{}}$\tilde{t}_{CPU}^{LOADS}$\\ \\\end{tabular}&
\begin{tabular}{@{}c@{}}$\varepsilon_{p}^{bd}$\\ (km)\\\end{tabular}&
\begin{tabular}{@{}c@{}}$\varepsilon_{v}^{bd}$\\ (m/s)\\\end{tabular}&
\begin{tabular}{@{}c@{}}$\varepsilon_{p}$\\ (km)\\\end{tabular}&
\begin{tabular}{@{}c@{}}$\varepsilon_{v}$\\ (m/s)\\\end{tabular}&
\begin{tabular}{@{}c@{}}$I_{\mu}$\\ ($\%$)\\\end{tabular}&
\begin{tabular}{@{}c@{}}$I_{\lambda}$\\ ($\%$)\\\end{tabular}&
\begin{tabular}{@{}c@{}}$I_{\gamma_{1,p}}$\\ ($\%$)\\\end{tabular}&
\begin{tabular}{@{}c@{}}$I_{\gamma_{2,p}}$\\ ($\%$)\\\end{tabular}\\
\noalign{\smallskip}\hline\noalign{\smallskip}
2& \phantom{0}349& 0.32& 3.55& 2.02& 5.35& 3.26& 2.31e-3& 5.02e-2& 1.038e1& 1.158e1\\
3& 1183& 1.00& 1.86& 1.12& 2.23& 1.24& 2.63e-4& 1.28e-2& 1.01e0& 1.49\phantom{0}\phantom{0}\phantom{0}\\
4& 2621& 1.60& 3.63& 2.27& 0.94& 0.53& 2.58e-4& 5.62e-3& 4.18e-3& 3.90e-2\\
\noalign{\smallskip}\hline
\end{tabular}
\end{center}
\end{table}
Starting from these considerations, the effect of the uncertainty set size on the described performance \hl{indices} was studied. The results are presented in Table~\ref{tab:LOADS_vs_MC_sigma}, for three values of $\alpha_{\lambda_j}$, namely 2, 3 (our reference solution) and 4. As expected, as the size of the uncertainty domain increases, the number of generated sets and the required computational load increase. Conversely, the accuracy of the statistical representation significantly improves, as shown by the \hl{indices} $I_{\mu}$, $I_\lambda$, $I_{\gamma_{1,p}}$ and $I_{\gamma_{2,p}}$. The most relevant improvement can be noticed in the last two \hl{indices} and can be explained by considering how the sample skewness and kurtosis are built. The two \hl{indices} compute a proper weighted sum of the deviation of the drawn samples from the mean of the distribution. As a result, the outliers are the elements with the largest contribution. When a $2\sigma$ representation is considered, most outliers fall out of the initial domain. Therefore, they are not included in any of the generated sets. To compute the \hl{indices}, they are assigned to the closest domain, but in any case they fall out of its [-1,1]$^n$ range. This implies a lower accuracy of the mapped samples, which generates values of $\gamma_{1,p}$ and $\gamma_{2,p}$ that are quite distant from the nominal \acrshort{mc} ones. When, instead, a larger initial domain is considered, the number of samples falling out of this domain decreases, thus the accuracy of the mapped samples increases, and $\gamma_{j,p}^{LOADS}\rightarrow\gamma_{j,p}^{MC}$. This is apparent if one compares the values of $I_{\gamma_{1,p}}$ and $I_{\gamma_{2,p}}$ for the $2\sigma$ and $4\sigma$ simulations, with a reduction in the relative errors of at least three orders of magnitude. Overall, the selection of the $\alpha_{\lambda_j}$ parameter shall be made according to the dimensionality of the problem as a trade-off between desired statistical accuracy and required computational load.

\subsection{\Acrshort{loads}-\acrshort{gmm}: test cases}
\label{subsec:LOADS-GMM_test}
This last section illustrates the application of the described \acrshort{loads}-\acrshort{gmm} algorithm to orbit uncertainty propagation. Two test cases taken from the literature are considered and used to assess the performance of our approach. The comparison is \hl{mostly} qualitative, since few absolute performance \hl{indices} are presented by the authors of these works.

The \acrshort{loads} algorithm is initialized by setting $N_{max,j}=10$ and $n_{b/r}=20$, while different values of $\varepsilon_{\nu}$ are analysed to assess the sensitivity on this parameter. The GMM pattern is built by considering the univariate splitting library described in Section~\ref{sec:LOADS-GMM}, with $\tilde{w}^{(0)}=\tilde{w}^{(2)}=0.2252246852539708$ and $\tilde{w}^{(1)}=0.5495506294920584$. The performance of the approach is expressed in terms of computational time and \acrfull{lam}~\citep{DeMars2013}. Given two multivariate \acrshortpl{pdf} $p$ and $q$, the \acrshort{lam} is expressed as
\begin{equation}
\mathcal{L}(p,q)=\int_{\mathbb{R}^n}p(\bm{x})q(\bm{x})\textrm{d}\bm{x}
\end{equation}
This index quantifies the overlap between the two distributions: the higher the \acrshort{lam}, the greater the agreement between them. In our analysis, we need to compare the estimated \acrshort{gmm} at the final epoch with the distribution resulting from the sampling of the initial Gaussian \acrshort{pdf} propagated pointwise to the final epoch. As a result, the \acrshort{lam} can be re-expressed as~\citep{DeMars2013}
\begin{equation}
\mathcal{L} = \dfrac{1}{N_{MC}}\sum\limits_{l=1}^{N_{MC}}\sum\limits_{k=1}^{n_{set}}w^{(k)} p _g\left(\bm{x}_{f,l};\bm{\mu}_f^{(k)},\bm{P}_f^{(k)}\right)
\end{equation}
with $N_{MC}$ number of \acrshort{mc} samples, $\bm{x}_{f,l}$ result of a pointwise propagation of the $l$-th \acrshort{mc} sample and $w^{(k)},\bm{\mu}_f^{(k)},\bm{P}_f^{(k)}$ the weights, means and covariances of the mixture model at $t_f$.

\begin{figure}[!t]
\centering
\subfloat[\label{subfig:DeMars_x-y}]{\includegraphics[trim=0cm 0cm 0cm 0cm, clip=true, width=0.7\textwidth]{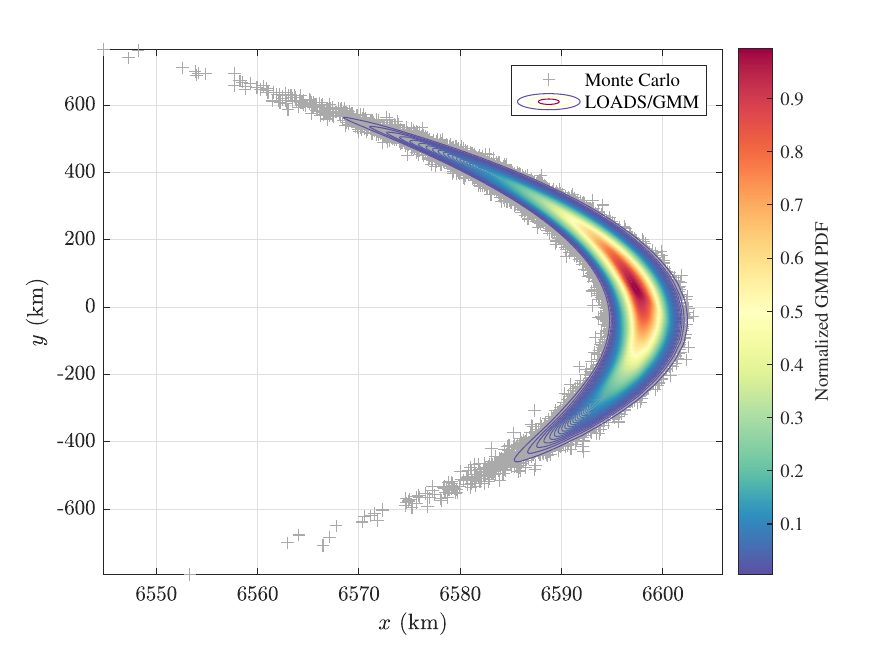}}\\
\subfloat[\label{subfig:DeMars_vx-vy}]{\includegraphics[trim=0cm 0cm 0cm 0cm, clip=true, width=0.7\textwidth]{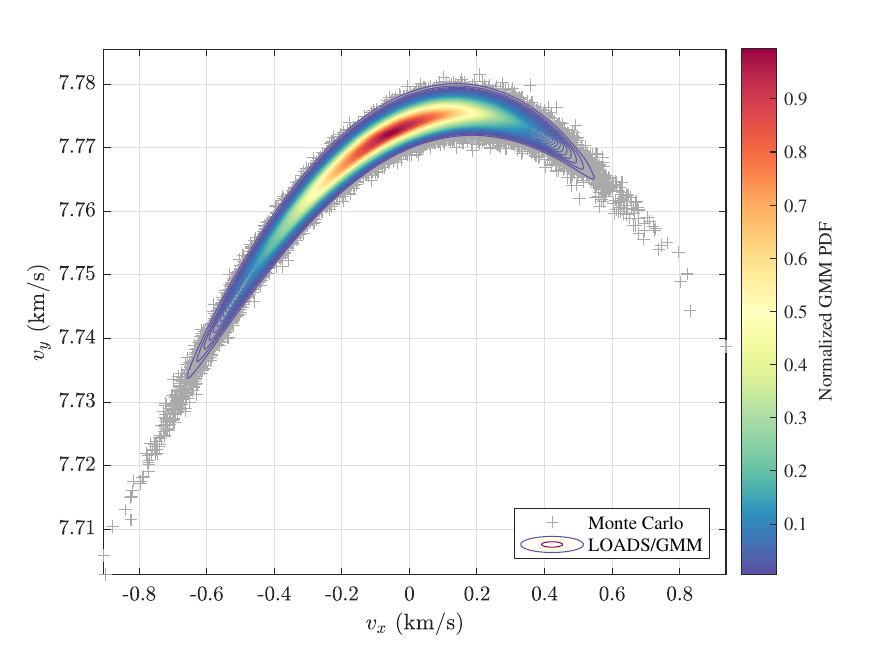}}\\
\caption{\Acrshort{mc} samples (grey crosses) and \acrshort{loads}-\acrshort{gmm} contours for the test case 1 after two orbital periods: \protect\subref{subfig:DeMars_x-y} $x$-$y$ plane and \protect\subref{subfig:DeMars_vx-vy} $v_x$-$v_y$ plane ($\varepsilon_{\nu}=0.03$, $N_{max,j}=10$, $n_{b/r}=20$, $L=3$, $\lambda = 1e-3$, 729 sets).}
\label{fig:LOADS-GMM_test1}
\end{figure}

The first test case is taken from~\citet{DeMars2013}. An object on a circular orbit \hl{at} 225~km of altitude is considered. The equations of motion are a simplified version of Eq.~\eqref{eq:EoM}, where the $J_2$ effect is neglected. The initial state distribution is assumed to be Gaussian, i.e. fully described with a covariance matrix, here considered diagonal, with standard deviations of $1.3$~km and $0.5$~km in $x$ and $y$ and $2.5$~m/s and $5$~m/s in $u$ and $v$. A propagation window of two orbital revolutions is considered. The AEGIS method proposed by~\citet{DeMars2013} adaptively splits the initial Gaussian distribution according to the nonlinearities encountered during the propagation. No merging algorithm is \hl{implemented}. As a result, \hl{around 300} sets are generated at the end of the propagation, while no clue on the required computational time is given. \hl{Moreover, the authors do not discuss how the splitting direction is chosen among the $n$ available ones.}

Figure~\ref{fig:LOADS-GMM_test1} shows the result of our \acrshort{loads}-\acrshort{gmm} method after a propagation of two orbital periods when considering $\varepsilon_{\nu}=0.03$. The algorithm generates 729 sets with a runtime of 3 seconds. The null eccentricity of the orbit prevents the \acrshort{loads} from performing any \hl{merging}, thus allowing for a fair comparison between our approach and the AEGIS method. The plots shows the comparison between \acrshort{mc} samples drawn from the initial Gaussian distribution and individually propagated to $t_f$ and the contours resulting from the propagated \acrshort{gmm} distribution $p_{\bm{X}_f}$. Two projections are considered, namely the $x-y$ plane (Fig.~\ref{subfig:DeMars_x-y}), and the $v_x-v_y$ plane (Fig.~\ref{subfig:DeMars_vx-vy}). As can be seen, the shape of the contours perfectly follows the curvature of the uncertainty region.

An analysis of the sensitivity of the algorithm on the $\varepsilon_{\nu}$ parameter is offered in Table~\ref{tab:LOADS-GMM_test1}. The table shows how the number of sets, \acrshort{lam} and computational time change as a function of the selected \acrshort{nli} threshold $\varepsilon_{\nu}$. Both \acrshort{lam} and runtime are normalized with respect to the $\varepsilon=0.03$ case. As for the \acrshort{loads} case, an increase in $\varepsilon_{\nu}$ causes a dramatic drop in the number of generated sets, and, consequently, in the required computational time. This however has an impact on the accuracy of the resulting \acrshort{gmm} distribution. A reduction of a factor 2 can be indeed noticed while passing from $\varepsilon_{\nu}=0.03$ to 0.04. As a result, the selection of the \acrshort{nli} threshold shall be a compromise between computational effort and desired accuracy. 

\hl{A direct comparison between the LOADS-GMM and the AEGIS method is not possible, since this should be done by fixing a reference accuracy (which is not available in absolute terms in~\citet{DeMars2013}) and comparing the number of generated sets, or the other way round. However, we can say that the obtained contours and required computational time look promising.}

\begin{table}[!t]
\caption{\Acrshort{loads}-\acrshort{gmm} performance for test case 1 as a function of $\varepsilon_{\nu}$ ($N_{max,j}=10$, $n_{b/r}=20$, $L=3$, $\lambda = 1e-3$).}
\label{tab:LOADS-GMM_test1}
\begin{center}
\begin{tabular}{*{4}{c}}
\hline\noalign{\smallskip}
$\varepsilon_{\nu}$& $n_{set}$& $\tilde{\mathcal{L}}$& $\tilde{t}_{CPU}$\\
\noalign{\smallskip}\hline\noalign{\smallskip}
3.0e-2& 729& 1.00& 1.00\\
3.5e-2& 243& 0.78& 0.49\\
4.0e-2& \phantom{0}93& 0.47& 0.30\\
$\infty$& \phantom{0}\phantom{0}1& 0.05& 0.08\\
\noalign{\smallskip}\hline
\end{tabular}
\end{center}
\end{table}

\begin{figure}[!t]
\centering
\subfloat[\label{subfig:Vish_x-y}]{\includegraphics[trim=0cm 0cm 0cm 0cm, clip=true, width=0.7\textwidth]{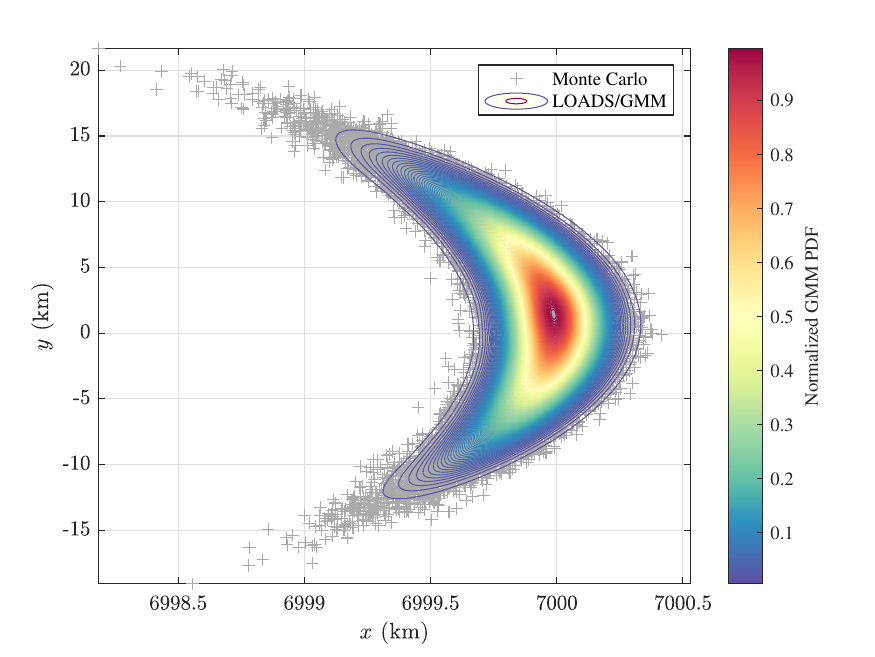}}\\
\subfloat[\label{subfig:Vish_vx-vz}]{\includegraphics[trim=0cm 0cm 0cm 0cm, clip=true, width=0.7\textwidth]{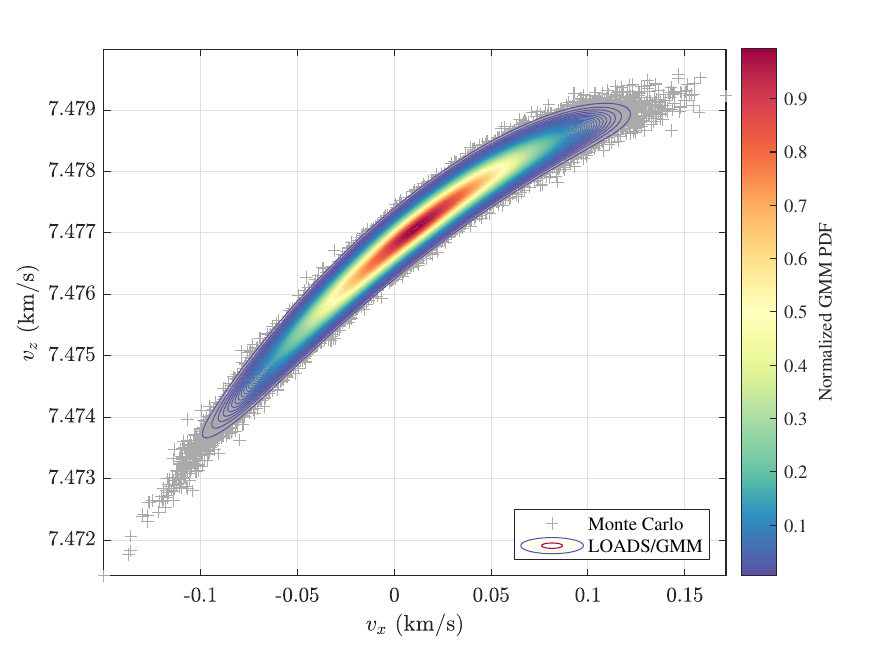}}\\
\caption{\Acrshort{mc} samples (grey crosses) and \acrshort{loads}-\acrshort{gmm} contours for the test case 2 after two orbital periods: \protect\subref{subfig:Vish_x-y} $x$-$y$ plane and \protect\subref{subfig:Vish_vx-vz} $v_x$-$v_z$ plane ($\varepsilon_{\nu}=0.01$, $N_{max,j}=10$, $n_{b/r}=20$, $L=3$, $\lambda = 1e-3$, 27 sets).}
\label{fig:LOADS-GMM_test2}
\end{figure}

\begin{table}[!t]
\caption{\Acrshort{loads}-\acrshort{gmm} performance for test case 2 as a function of $\varepsilon_{\nu}$ ($N_{max,j}=10$, $n_{b/r}=20$, $L=3$, $\lambda = 1e-3$).}
\label{tab:LOADS-GMM_test2}
\begin{center}
\begin{tabular}{*{4}{c}}
\hline\noalign{\smallskip}
$\varepsilon_{\nu}$& $n_{set}$& $\tilde{\mathcal{L}}$& $\tilde{t}_{CPU}$\\
\noalign{\smallskip}\hline\noalign{\smallskip}
7.0e-3& 243& 3.15& 3.01\\
1.0e-2& \phantom{0}27& 1.00& 1.00\\
1.5e-2& \phantom{0}\phantom{0}3& 0.21& 0.64\\
$\infty$& \phantom{0}\phantom{0}1& 0.10& 0.56\\
\noalign{\smallskip}\hline
\end{tabular}
\end{center}
\end{table}

A second test case is taken from~\citet{Vishwajeet2018}. In this case the authors implement a split/merge technique. The test case considers the motion of a Low Earth Orbit satellite subject to Keplerian motion only. At the initial epoch, the states are assumed to have a Gaussian distribution with the following mean $\bm{\mu}_0$ and covariance $\bm{P}_0$
\begin{equation}
\begin{gathered}
\bm{\mu}_0 = \left\{7000,0,0,0,-1.0374,7.4771\right\}^{\textrm{T}}\\
\bm{P}_0=\textrm{diag}\left\{0.01,0.01,0.01,1\textrm{e}-6,1\textrm{e}-6,1\textrm{e}-6\right\}
\end{gathered}
\end{equation}
where km and km/s are considered for position and velocity, respectively. The state is then propagated for 3.25~h. The authors obtain around 150 sets in 2~h of runtime.

Figure~\ref{fig:LOADS-GMM_test2} shows the result of our \acrshort{loads}-\acrshort{gmm} method after 3.25~h of propagation when considering $\varepsilon_{\nu}=0.01$. The algorithm generates 27 sets in just 491 milliseconds. The plots show the comparison between \acrshort{mc} samples and \acrshort{gmm} contours on the $x-y$ (Fig.~\ref{subfig:Vish_x-y}) and $v_x-v_z$ (Fig.~\ref{subfig:Vish_vx-vz}) planes. \hl{Similarly as before}, a very good matching between the curvature of the contours and the shape of the propagated uncertainty region is obtained. Table~\ref{tab:LOADS-GMM_test2} finally shows the sensitivity of the method on $\varepsilon_{\nu}$. In this case, a stronger variability of the \acrshort{lam} is obtained. 

\hl{Even though a rigorous comparison between the two approaches is not possible, the extremely low computational burden and very good matching between \acrshort{gmm} contours and \acrshort{mc} distribution corroborate the expectations from the previous simulation and lead us to conclude that our approach shows appealing performance and could be considered a valuable tool for uncertainty propagation}.

\section{Conclusions}
This paper introduced a \acrfull{loads} algorithm for the nonlinear transformation of sets. The algorithm is built on a \acrfull{nli} that exploits \acrfull{da} and polynomial bounding techniques to estimate the nonlinear nature of a problem. More specifically, starting from a generic nonlinear transformation and uncertainty set described with a second-order Taylor expansion, \acrshort{da} provides the first-order Taylor expansion of the terms of the Jacobian matrix of the transformation. By exploiting polynomial bounding techniques, the upper and lower bounds of each term can be rigorously computed, and they can be used to build an estimate of the nonlinearity of the transformation. The proposed formulation removes the need for sampling, thus helping the user investigate the uncertainty set and avoiding the risk of inaccurate estimates.  The described index is then embedded into an algorithm that accurately maps an uncertainty set through a nonlinear transformation by progressively splitting and approximating it as a manifold of second-order Taylor expansions. 

The first part of the paper describes the algorithm, including its tailoring to the case of orbital uncertainty propagation. The \acrshort{loads} algorithm is coupled with a merging technique that minimizes the number of propagated sets according to the current level of nonlinearities. A coupling of the \acrshort{loads} algorithm with \hl{a} \acrfull{gmm} \hl{approximation} of the \acrshort{pdf} is then proposed as a valuable tool to obtain an analytical description of the propagated state \acrshort{pdf}. In the second part of the paper the performance of \acrshort{nli}, \acrshort{loads} and \acrshort{loads}-GMM algorithms is demonstrated in the context of orbit uncertainty propagation. The analyses include an extensive investigation of the role of the control parameters on the accuracy and efficiency of the \acrshort{loads} algorithm. The \acrshort{loads} algorithm is an effective approach for the nonlinear propagation of uncertainty with an accuracy and computational load that can be controlled by tuning the nonlinearity threshold and the \hl{size of the} uncertainty domain. The availability of local low-order maps represents one of the strongest points of \hl{this} approach, as they can be used to map \hl{the} statistics with techniques that would be inaccurate on the whole uncertainty set.  The coupling of the \acrshort{loads} algorithm with \hl{\acrshortpl{gmm}} represents a natural extension of the approach. The \acrshort{loads} algorithms \hl{provides in fact the tools to adaptively increase the number of Gaussian kernels online and recombine them} whenever possible. The results show that \hl{this} method has \hl{appealing} performance, both in accuracy and required computational time.

Future work will be dedicated to the application of the developed \hl{\acrshort{loads}-GMM method} to the propagation of large uncertainty sets resulting from short-arc initial orbit determination, and the combination with coordinates sets that limit the nonlinearity in orbital dynamics.


\section*{Acknowledgements}

This work made use of the CNES orbital propagation tools, including the PACE library.

\bibliography{LOADS.bib}

\end{document}